\documentclass[12pt]{article}
\usepackage{amsthm, amsmath,amssymb,amsfonts}
\usepackage{graphicx}
\usepackage{color}

\usepackage{psfig}
\usepackage{graphics}
\usepackage{amsmath}
\usepackage{amssymb}
\usepackage{epsfig}
\usepackage{epstopdf}
\usepackage{color}

\newtheorem{theorem}{Theorem}[section]

\newtheorem{lemma}[theorem]{Lemma}
\newtheorem{remark}[theorem]{Remark}
\newtheorem{definition}[theorem]{Definition}
\newtheorem{example}[theorem]{Example}

\def\para{\vspace{2mm}}
\def\cP{{\mathcal P}}

\def\gcd{{\rm gcd}}
\def\lcm{{\rm lcm}}

\def\ord{{\rm ord}}

\def\deg{{\rm deg}}

\def\limit{{\rm limit}}
\def\ox{\,{\overline x}\,}
\def\Re{{\rm Re}}
\def\Im{{\rm Im}}

\begin{document}
\title{Characterizing the finiteness of the Hausdorff distance between two algebraic curves\thanks{The author S. P\'erez-D\'{\i}az is member of the Research Group ASYNACS (Ref. CCEE2011/R34)}}
\author{Angel Blasco and Sonia P\'erez-D\'{\i}az\\
Departamento de F\'{\i}sica y Matem\'aticas \\
        Universidad de Alcal\'a \\
      E-28871 Madrid, Spain  \\
angel.blasco@uah.es, sonia.perez@uah.es
}
\date{}          
\maketitle

\begin{abstract}
In this paper, we present a characterization for the Hausdorff
distance between two given algebraic curves in the $n$-dimensional
space (parametrically or implicitly defined) to be finite. The
characterization is related with the asymptotic behavior of the two
curves and it can be easily checked. More precisely,  the Hausdorff
distance between two curves $\cal C$ and $\overline{\cal C}$ is
finite if and only if for each infinity branch of $\cal C$ there exists an
infinity branch of $\overline{\cal C}$ such that the terms with
positive exponent in the corresponding series are the same, and
reciprocally.


\end{abstract}

{\bf Keywords:} Hausdorff Distance; Algebraic Space Curves; Implicit
Polynomial; Parametrization; Infinity Branches; Asymptotic Behavior;


 \section{Introduction}

The Hausdorff distance is one of the most used measures in geometric pattern matching algorithms, computer aided design or computer graphics (see e.g. \cite{Hutt}, \cite{Kim}, \cite{Patri}, \cite{Vivek}).

\para

  Intuitively speaking, given a metric space $(E,d)$ and two  arbitrary subsets $A, B\subset E$, the  Hausdorff distance
assigns to each point of one set the distance to its closest point on the other and takes
the maximum over all these values  (see \cite{Ali}). More precisely, the   {\em Hausdorff distance} between $A$ and $B$ is defined
as:
$$d_H(A,B)=\max\{\sup_{x\in A}\inf_{y\in B}d(x,y),\sup_{y\in B}\inf_{x\in A}d(x,y)\}.$$
In this paper, we deal with the particular case where  $E=\mathbb{C}^n$, $d$ is the usual unitary   distance, and the two arbitrary subsets
are  two real algebraic  curves ${\cal C}$ and
${\overline{{\cal C}}}$. In this case, the Hausdorff distance between ${\cal C}$ and
${\overline{{\cal C}}}$ is given by
$$d_H({\cal C},{\overline{{\cal C}}})=\max\{\sup_{p\in {\cal C}}d(p,{\overline{{\cal C}}}),\sup_{\overline{p}\in {\overline{{\cal C}}}}d(\overline{p},{\cal C})\}$$
where $d(p,{\cal C})=\min\{d(p,q):q\in{\cal C}\}$.\\

In general,  $d_H(A,B)$ may be infinite, and some restrictions have
to be  imposed  to guarantee its finiteness (see e.g. \cite{RSS}).

\para

As far as the authors know, there is no efficient algorithms for the
exact computation of the Hausdorff distance between algebraic
varieties (in fact, if both varieties are given in implicit form,
the computation of the Hausdorff distance is even harder). Only some
results for bounding or estimating the Hausdorff distance as well as
computing it for some special cases can be found (see e.g.
\cite{Bai}, \cite{Chen}, \cite{Henri},    \cite{Juttler},
\cite{Kim}, \cite{RSS2}). These results  play an important role in
some applications to computer aided geometric design as for instance
in  the approximate parametrization problem (see e.g. \cite{PSS1},
\cite{PSS2},  \cite{PRSS}, \cite{RS}, \cite{RSS}).   In that
problem,   given an affine curve $\cal C$ (say that it is a
perturbation of a rational curve), the goal is to compute a rational
parametrization of a rational affine curve $\overline{{\cal C}}$
near $\cal C$ (one may state the problem also for surfaces). The
effectiveness of the algorithm will depend on the closeness of
${\cal C}$ and $\overline{{\cal C}}$ and, at least, one needs to
show that  the Hausdorff distance between ${\cal C}$ and
$\overline{{\cal C}}$ is finite. The potential applications of the
Hausdorff distance also include error bounds for the approximate
implicitization of  curves and surfaces (see e.g. \cite{Bajaj},
\cite{Doken}, \cite{Emiris}).


\para

In this paper, we  characterize whether the Hausdorff distance
between two given algebraic curves in the $n$-dimensional space is
finite. These two curves can be both, parametrically or implicitly
defined. The characterization improves Proposition 5.4 in
\cite{paper1}, and it is based on the notion of infinity branch
which reflects the status of a curve at the points with sufficiently
large coordinates.

\para

This concept is   an essential  tool to analyze the
 behavior at the infinity of an   algebraic   curve, which implies a wide applicability in many active research fields. For instance, infinity branches  allow us  to  sketch  the graph of a given algebraic curve as well as to study  its topology (see e.g. \cite{Gao}, \cite{lalo}, \cite{Hong}). In addition, the notion of g-asymptote is introduced from the concept of infinity branch (see \cite{BlascoPerezII} and \cite{paper3}). We say that  a curve $\overline{{\cal C}}$ is a {\it  generalized asymptote} (or \textit{g-asymptote}) of another curve $\cal C$ if
$\overline{{\cal C}}$ approaches $\cal C$ at some infinity branch, and $\cal C$  can not be approached at that branch by a new curve of lower degree (that is, the notion of  g-asymptote generalizes  the classical notion of (linear) asymptote).

\para

In this paper, we introduce the concept of curves, $\cal C$ and
$\overline{\cal C}$,  having a \textit{similar asymptotic behavior},
which is concerned with the convergence/divergence of their infinity
branches. More precisely, we say that   $\cal C$ and $\overline{\cal
C}$ have a \textit{similar asymptotic behavior} if there are no
infinity branches in $\cal C$ which diverge from all the infinity
branches in $\overline{\cal C}$, and reciprocally.



\para

From this concept, we prove the main theorem, which states a
necessary and sufficient condition for the Hausdorff distance
between two curves to be finite.
More precisely, we show that, given two algebraic curves in the
affine $n$-space, the Hausdorff distance between them is finite if
and only if they have {\it a similar asymptotic behavior}. This
condition is very easy to formulate from the  computational point of view and thus, we
present an effective algorithm that checks if it holds.

\para

The structure of the paper is as follows: In Section 2, we present
the terminology that will be used throughout the paper as well as
some previous results. These results are presented for both, curves
given implicitly and curves defined parametrically. Section 3, is
devoted to present the main theorem where the finiteness of the
Hausdorff distance is characterized. For this purpose, some previous
technical lemmas are proved. In addition, we derive an algorithm
that determine whether the Hausdorff distance between two given
algebraic curves is finite and we illustrate it with some examples.

 \section{Notation and terminology}\label{S-notation}

In this section, we present some notions and terminology that will
be used throughout the paper. In particular, we need some previous
results concerning local parametrizations and Puiseux series.  For
further details see \cite{Alon92},  \cite{paper1}, \cite{Duval89},
Section 2.5 in \cite{SWP},    and Chapter 4 (Section 2)
in \cite{Walker}.

\para

We denote by ${\Bbb C}[[t]]$ the domain of {\em formal power series}
in the indeterminate $t$ with coefficients in the field ${\Bbb C}$,
i.e. the set of all sums of the form $\sum_{i=0}^\infty a_it^i$,
 $a_i \in {\Bbb C}$. The quotient field of ${\Bbb C}[[t]]$ is
called the field of {\em formal Laurent series}, and it is denoted by
${\Bbb C}((t))$. It is well known that every non-zero formal Laurent
series $A \in {\Bbb C}((t))$ can be written in the form
$A(t) = t^k\cdot(a_0+a_1t+a_2t^2+\cdots), {\rm\ where\ }
    a_0\not= 0 {\rm\ and\ } k\in \Bbb Z.$
  In addition, the field
${{\Bbb C}\ll t\gg}:= \bigcup_{n=1}^\infty {\Bbb C}((t^{1/n}))$ is
called the field of {\em  formal Puiseux series}. Note that Puiseux
series are power series of the form
$$\varphi(t)=m+a_1t^{N_1/N}+a_2t^{N_2/N}+a_3t^{N_3/N} +
\cdots\in{\Bbb C}\ll t\gg,\quad a_i\not=0,\, \forall i\in {\Bbb
N},$$
 where  $N, N_i\in {\Bbb N},\,\,i\geq 1$,  and $0<N_1<N_2<\cdots$. The natural number $N$ is known as
 {\em  the ramification index} of the series. We denote it as $\nu(\varphi)$ (see \cite{Duval89}).

\para

The {\em order} of a non-zero (Puiseux or Laurent) series $\varphi$ is the
smallest exponent of a term with non-vanishing coefficient in $\varphi$.
We denote  it  by $\ord(\varphi)$. We let the order of 0 be $\infty$.

\para

The most important property of Puiseux series is given by Puiseux's Theorem, which states that if $\Bbb K$ is an algebraically closed field, then the field ${\Bbb K}\ll x\gg$ is algebraically closed  (see Theorems 2.77 and 2.78 in \cite{SWP}). A proof of Puiseux's Theorem can be given constructively by the
Newton Polygon Method (see e.g. Section 2.5 in \cite{SWP}).

\para

In the following, we deal with space curves that are  implicitly defined. In Subsection \ref{Sub-Param}, we will consider space curves parametrically defined.

\subsection{Implicitly defined space curves} \label{Sub-Impli}

Let ${\cal C}\in {\Bbb C}^{n}$ be a curve in the $n$-dimensional space defined
by a finite set of real polynomials $f_1(\ox),\ldots, f_s(\ox)\in {\Bbb R}[\ox],\,s\geq n-1$, where $\ox=(x_1,\ldots,x_n)$.

\para

The assumption of reality of the curve $\cal C$ is included because of the nature of the
problem, but the theory developed in this paper can be applied  for the case of
complex non-real curves.

\para

 Let
${\cal C}^*$  be the corresponding projective curve  defined by the
homogeneous polynomials $F_i(x_1,\ldots,x_n,x_{n+1}) \in {\Bbb
R}[x_1,\ldots,x_n,x_{n+1}],\,i=1,\ldots,s$. Furthermore, let $P=(1:m_2:\ldots: m_n:0),\,m_j\in {\Bbb
C},\,j=2,\ldots,n$ be an infinity point of ${\cal C}^*$.

\para

In addition,  we consider the curve
 implicitly defined by the polynomials $g_i(x_2,\ldots,x_n,x_{n+1}):=F_i(1,x_2,\ldots,x_n,x_{n+1})\in
{\Bbb R}[x_2,\ldots,x_n,x_{n+1}]$ for $i=1,\ldots,s$. Observe that
$g_i(p)=0,$ where $p=(m_2,\ldots, m_n,0)$.  Let $I\in {\Bbb
R}(x_{n+1})[x_2,\ldots,x_n]$ be the ideal generated by
$g_i(x_2,\ldots,x_n,x_{n+1}),\,\,i=1,\ldots,s,$ in the ring ${\Bbb
R}(x_{n+1})[x_2,\ldots,x_n]$. We assume that ${\cal C}$ is not
contained in some hyperplane $x_{n+1} =c,\,c\in {\Bbb C}$ (otherwise, one can consider ${\cal C}$ as a curve in the $(n-1)$-dimensional space), and thus
we have that  $x_{n+1}$ is not algebraic over $\Bbb R$. Under this
assumption, the ideal $I$ (i.e. the system of equations $g_1 =
\cdots = g_s = 0$) has only finitely many solutions in the
$n$-dimensional affine space over the algebraic  closure of ${\Bbb
R}(x_{n+1})$ (which is contained in ${{\Bbb C}\ll x_{n+1}\gg}$).
Then, there are finitely many $(n-1)$-tuples
$(\varphi_{2}(t),\ldots,\varphi_{n}(t))$ where $\varphi_{j}(t)\in
{{\Bbb C}\ll t\gg},\,j\in \{2,\ldots,n\}$, such that
$g_i(\varphi_{2}(t),\ldots,\varphi_{n}(t),t) = 0,\,i=1,\ldots,s$,
 and $\varphi_{j}(0)=m_j,\,j=2,\ldots,n$. Each of these
$(n-1)$-tuples is a solution of the system associated
with the infinity point $(1:m_2:\ldots:m_n:0)$,  and each
$\varphi_{j}(t)$ converges in a neighborhood of $t= 0$. Moreover,
 since $\varphi_{j}(0)=m_j,\,j=2,\ldots,n$, these series
do not have terms with negative exponents; in fact, they have the
form
$$\varphi_{j}(t)=m_j+\sum_{i\geq 1}a_{i,j} t^{N_{i,j}/N_j}$$
where $N_j,\,N_{i,j}\in {\Bbb N},\,\,0<N_{1,j}<N_{2,j}<\cdots,$.

\para

It is important to remark that if
$\varphi(t):=(\varphi_{2}(t),\ldots,\varphi_{n}(t))$ is a solution
of the system, then
$\sigma_{\epsilon}(\varphi)(t):=(\sigma_{\epsilon}(\varphi_{2})(t),\ldots,\sigma_{\epsilon}(\varphi_{n})(t))$
is another solution of the system, where
$$\sigma_{\epsilon}(\varphi_{j})(t)=m_j+\sum_{i\geq 1}a_{i,j}\epsilon^{\lambda_{i,j}} t^{N_{i,j}/N_j},\,\,N_j,\,N_{i,j}\in {\Bbb N},\,\,0<N_{1,j}<N_{2,j}<\cdots,$$
$N:=\lcm(N_2,\ldots, N_n)$, $\lambda_{i,j}:=N_{i,j}N/N_j\in {\Bbb
N}$, and $\epsilon^N=1$ (see \cite{Alon92}). We refer to these
solutions as the {\em conjugates} of $\varphi$. The set of all
(distinct) conjugates of $\varphi$ is called the  {\em conjugacy
class} of $\varphi$, and the number of different conjugates
is $N$.  We denote the natural number $N$  as  $\nu(\varphi)$.

\para

Under these conditions and reasoning as in \cite{paper1}, we get
that there exists   $M \in {\Bbb R}^+$ such that for $i=1,\ldots,s$,
$$F_i(1:\varphi_{2}(t):\ldots:\varphi_{n}(t):t)=g_i(\varphi_{2}(t),\ldots,\varphi_{n}(t),t)=0$$
for\, $t\in {\Bbb C}$\, and $|t|<M$.
This implies that
$$F_i(t^{-1}:t^{-1}\varphi_{2}(t):\ldots: t^{-1}\varphi_{n}(t):1)=f_i(t^{-1},t^{-1}\varphi_{2}(t),\ldots, t^{-1}\varphi_{n}(t))=0,$$ for
$t\in {\Bbb C}$  and  $0<|t|<M$.

\para

\noindent
Now, we set $t^{-1}=z$, and   we obtain
that  for $i=1,\ldots,s$,
$$f_i(z,r_{2}(z),\ldots,r_{n}(z))=0,\quad \mbox{$z\in {\Bbb C}$\, and
$|z|>M^{-1}$,\qquad where}$$
\[r_{j}(z)=z\varphi_{j}(z^{-1})=\]
\begin{equation}\label{Eq-r_j}
m_jz+a_{1,j}z^{1-N_{1,j}/N_j}+a_{2,j}z^{1-N_{2,j}/N_j}+a_{3,j}z^{1-N_{3,j}/N_j}
+ \cdots,\end{equation} $a_{i,j}\not=0$, $N_j,N_{i,j}\in {\Bbb
N},\,\,i=1,\ldots$, and $0<N_{1,j}<N_{2,j}<\cdots$.

\para

\para

  Since $\nu(\varphi)=N$, we get that there are $N$ different series in its conjugacy
class. Let  $\varphi_{\alpha,j},\,\alpha=1,\ldots,N$ be these
series, and $$r_{\alpha,j}(z)=z\varphi_{\alpha,j}(z^{-1})=$$
\begin{equation}\label{Eq-conjugates}
m_jz+a_{1,j}c_\alpha^{\lambda_{1,j}}
z^{1-N_{1,j}/N_j}+a_{2,j}c_\alpha^{\lambda_{2,j}}
z^{1-N_{2,j}/N_j}+a_{3,j}c_\alpha^{\lambda_{3,j}} z^{1-N_{3,j}/N_j}
+ \cdots
\end{equation}
where $N:=\lcm(N_2,\ldots,N_n)$, $\lambda_{i,j}:=N_{i,j}N/N_j\in
{\Bbb N}$, and $c_1,\ldots,c_N$ are the $N$ complex roots of
$x^N=1$. Now we are ready to introduce the notion of infinity
branch. The following definitions and results generalize those
presented in \cite{paper1} for algebraic plane curves, and in
\cite{paper3} for algebraic space curves.

\para

\begin{definition}\label{D-infinitybranch}
 An {\em
infinity branch of a $n$-dimensional space  curve ${\cal C}$}
associated to the infinity point $P=(1:m_2:\ldots: m_n:0),\,m_j\in
{\Bbb C},\,j=2,\ldots,n$, is  a set $ B=\bigcup_{\alpha=1}^N
L_\alpha$, where $L_\alpha=\{(z,r_{\alpha,2}(z),\ldots,
r_{\alpha,n}(z))\in {\Bbb C}^n: \,z\in {\Bbb C},\,|z|>M\}$,\,  $M\in
{\Bbb R}^+$, and the series $r_{\alpha,j},\,j=2,\ldots,n,$ are given
by (\ref{Eq-conjugates}).
 The subsets $L_1,\ldots,L_N$ are
called the {\em leaves} of the infinity branch $B$.
 \end{definition}

\para

\begin{remark} \label{R-conjugation}
An infinity branch  is uniquely determined  from one leaf, up to
conjugation. That is, let $B$ be an infinity branch and let
$$L=\{(z,r_{2}(z),\ldots, r_{n}(z))\in {\Bbb C}^n: \,z\in {\Bbb
C},\,|z|>M\}$$ be one of its leaves, with
$$r_{j}(z)=z\varphi_{j}(z^{-1})=
m_jz+a_{1,j} z^{1-N_{1,j}/N_j}+a_{2,j} z^{1-N_{2,j}/N_j}+a_{3,j}
z^{1-N_{3,j}/N_j} + \cdots.$$ Then, any other leaf $L_\alpha$ has
the form
$$L_\alpha=\{(z,r_{\alpha,2}(z),\ldots, r_{\alpha,n}(z))\in {\Bbb
C}^n: \,z\in {\Bbb C},\,|z|>M\}$$ where
$r_{\alpha,j}=r_{j},\,j=2,\ldots,N$, up to conjugation; i.e.
$$r_{\alpha,j}(z)=z\varphi_{\alpha,j}(z^{-1})=$$$$=m_jz+a_{1,j}c_\alpha^{\lambda_{1,j}}
z^{1-N_{1,j}/N_j}+a_{2,j}c_\alpha^{\lambda_{2,j}}
z^{1-N_{2,j}/N_j}+a_{3,j}c_\alpha^{\lambda_{3,j}} z^{1-N_{3,j}/N_j}
+ \cdots$$ $N, N_{i,j}\in\mathbb{N}$,
$\lambda_{i,j}:=N_{i,j}N/N_j\in {\Bbb N},\,j=2,\ldots,n$ and
$c_{\alpha}^N=1,\,\,\alpha=1,\ldots,N$.
\end{remark}

\para

\begin{remark} \label{R-infinitypoint} Observe that the above approach is presented for infinity points of the
form $(1: m_2:\ldots: m_n: 0)$. For the infinity points $(0:
m_2:\ldots: m_n: 0)$, with $m_j\not=0$ for some $j=2,\ldots, n$, we
reason similarly but we dehomogenize w.r.t $x_j$. More precisely,
let us assume that $m_2\not=0$. Then,  we consider  the curve
defined by the polynomials $g_i(x_1,x_3,\ldots,
x_{n+1}):=F_i(x_1,1,x_3,\ldots, x_{n+1})\in {\Bbb R}[x_1,x_3,\ldots,
x_{n+1}],\,i=1,\ldots,s$, and we reason as above. We get that   an
infinity branch of   ${\cal C}$ associated to the infinity point
$P=(0:m_2:\ldots: m_n:0),\,m_2\not=0$, is  a set $
B=\bigcup_{\alpha=1}^N L_\alpha$, where
$L_\alpha=\{(r_{\alpha,1}(z),z,r_{\alpha,3}(z),\ldots,r_{\alpha,n}(z))\in
{\Bbb C}^n: \,z\in {\Bbb
C},\,|z|>M\}$,\,  $M\in {\Bbb R}^+$. \\
Additionally, instead of working with this type of branches, if the  space curve $\cal C$
 has infinity points of the form $(0 : m_2 : \ldots: m_n : 0)$, one may consider a linear change of coordinates. Thus, in the following, we may  assume w.l.o.g that the given algebraic  curve $\cal C$ only  has infinity points of the form $(1 : m_2 : \ldots: m_n : 0)$. More details on these type of branches are given in \cite{paper1} and \cite{paper3}.

\end{remark}

\para

In the following, we introduce the notions of convergent and
divergent leaves. Intuitively speaking, two leaves converge
(diverge) if they get closer (get away) as they tend to infinity.

\begin{definition}\label{D-distance0}
Let $L=\{(z,r_{2}(z),\ldots, r_{n}(z))\in {\Bbb C}^n:\,z\in {\Bbb
C},\,|z|>M\}$ and $\overline{L}=\{(z,\overline r_{2}(z),\ldots,
\overline r_{n}(z))\in {\Bbb C}^n:\,z\in {\Bbb
C},\,|z|>\overline{M}\}$ be two leaves that belong to two infinity
branches $B$ and $\overline{B}$, respectively. We say that
\begin{enumerate}
\item  $L$ and
$\overline{L}$ {\sf converge} if
$$\lim_{z\rightarrow\infty} d(({r}_2(z),\ldots,{r}_n(z)),
(\overline{r}_2(z),\ldots,\overline{r}_n(z)))=0.$$
\item $L$ and
$\overline{L}$ {\sf diverge} if
$$\lim_{z\rightarrow\infty} d(({r}_2(z),\ldots,{r}_n(z)),
(\overline{r}_2(z),\ldots,\overline{r}_n(z)))=\infty.$$
\end{enumerate}
\end{definition}

\begin{remark} \label{R-distance0}
We consider any distance $d(u,v)=\|u-v\|,\,u,v\in {\Bbb C}^{n-1}$,
where $\|p\|$ denotes the norm of a point $p\in  {\Bbb C}^{n-1}$. We
recall that all norms are equivalent in ${\Bbb C}^{n-1}$. Hence,
\begin{enumerate}
\item $\lim_{z\rightarrow\infty} d(({r}_2(z),\ldots,{r}_n(z)),
(\overline{r}_2(z),\ldots,\overline{r}_n(z)))=0$ if and only if\\
$\lim_{z\rightarrow\infty} (\overline{r}_j(z)-r_j(z))=0$ for every
$j=2,\ldots,n$.
\item $\lim_{z\rightarrow\infty}
d(({r}_2(z),\ldots,{r}_n(z)),
(\overline{r}_2(z),\ldots,\overline{r}_n(z)))=\infty$ if and only
if\\ $\lim_{z\rightarrow\infty} (\overline{r}_j(z)-r_j(z))=\infty$
for some $j=2,\ldots,n$.
\end{enumerate}
\end{remark}

\para

\begin{remark}\label{R-limit-c}
Observe that it may happen that
$$\lim_{z\rightarrow\infty} d(({r}_2(z),\ldots,{r}_n(z)),
(\overline{r}_2(z),\ldots,\overline{r}_n(z)))=c\in\mathbb{R}^+\setminus \{0\} $$
which is equivalent to  $\lim_{z\rightarrow\infty}
(\overline{r}_j(z)-r_j(z))=c_j\in {\Bbb C}$ for every $j=2,\ldots,n$
and $c_j\neq 0$ for some $j=2,\ldots,n$. In this case, $L$ and
$\overline{L}$   {\sf do not converge neither diverge} (compare with Definition \ref{D-distance0}).
\end{remark}

\para

The following lemma provides a procedure to determine whether two leaves
converge or diverge without the need of computing limits.

\begin{lemma}\label{L-DistVertical}
Let $L=\{(z,r_{2}(z),\ldots, r_{n}(z))\in {\Bbb C}^n:\,z\in {\Bbb
C},\,|z|>M\}$ and $\overline{L}=\{(z,\overline r_{2}(z),\ldots,
\overline r_{n}(z))\in {\Bbb C}^n:\,z\in {\Bbb
C},\,|z|>\overline{M}\}$ be two leaves that belong to two infinity
branches $B$ and $\overline{B}$, respectively. It holds that:
\begin{enumerate}
\item  $L$ and $\overline{L}$ converge if and only if the terms with
non-negative exponent in the series $r_j(z)$ and $\overline{r}_j(z)$
are the same, for every $j=2,\ldots, n$.
\item $L$ and $\overline{L}$ diverge if and only if the terms with
positive exponent in the series $r_j(z)$ and $\overline{r}_j(z)$ are
not the same, for some $j=2,\ldots, n$.
\end{enumerate}
\end{lemma}

\noindent\textbf{Proof:} Let
$$r_{j}(z)=
m_jz+a_{1,j} z^{1-N_{1,j}/N_j}+a_{2,j} z^{1-N_{2,j}/N_j}+a_{3,j}
z^{1-N_{3,j}/N_j} + \cdots,$$ $a_{i,j}\not=0,\, \forall i\in {\Bbb
N},\,i\geq 1,$ $N_j, N_{i,j}\in {\Bbb N}$, and
$0<N_{1,j}<N_{2,j}<\cdots$ for $j=2,\ldots, n$. and
$$\overline{r}_{j}(z)=
\overline{m}_jz+\overline{a}_{1,j}
z^{1-\overline{N}_{1,j}/\overline{N}_j}+\overline{a}_{2,j}
z^{1-\overline{N}_{2,j}/\overline{N}_j}+\overline{a}_{3,j}
z^{1-\overline{N}_{3,j}/\overline{N}_j} + \cdots,$$
$\overline{a}_{i,j}\not=0,\, \forall i\in {\Bbb N},\,i\geq 1,$
$\overline{N}_j, \overline{N}_{i,j}\in {\Bbb N}$, and
$0<\overline{N}_{1,j}<\overline{N}_{2,j}<\cdots$ for $j=2,\ldots,
n$. Then,
$$r_j(z)-\overline{r}_j(z)=m_jz-\overline{m}_jz+a_{1,j}z^{\frac{N-N_1}{N}}-\overline{a}_{1,j}z^{\frac{\overline{N}-\overline{N}_1}{\overline{N}}}
+a_{2,j}z^{\frac{N-N_2}{N}}-\overline{a}_{2,j}z^{\frac{\overline{N}-\overline{N}_2}{\overline{N}}}+\cdots.$$
Under these conditions, it holds that:
\begin{enumerate}
\item $\lim_{z\rightarrow\infty}(r_j(z)-\overline{r}_j(z))=0$ for
every $j=2,\ldots,n$,  if and only if all the exponents in the series
$r_j(z)-\overline{r}_j(z)$ are negative. This situation holds if the
terms with non-negative exponent in the series  $r_j(z)$ and
$\overline{r}_j(z)$  are the same for every $j=2,\ldots,n$.
\item $\lim_{z\rightarrow\infty}(r_j(z)-\overline{r}_j(z))=\infty$
for some $j=2,\ldots,n$, if and only if $r_j(z)-\overline{r}_j(z)$
has some term with positive exponent. This situation holds if the
terms with positive exponent in the series, $r_j(z)$ and
$\overline{r}_j(z)$, are not the same for some $j=2,\ldots,n$.\hfill $\Box$
\end{enumerate}

\begin{remark}
If the terms with positive exponent in the series $r_j(z)$ and
$\overline{r}_j(z)$ are the same for every $j=2,\ldots, n$, but the
independent terms (the terms with exponent zero) are different for some
$j=2,\ldots, n$, we have that $L$ and $\overline{L}$  do not diverge neither
converge.
\end{remark}

In the following, we introduce the notions of convergent and
divergent branches. These concepts are obtained from Definition \ref{D-distance0}, and they are an indispensable tool for comparing the asymptotic
behavior of two curves.

\para

\begin{definition}\label{D-convergent-branches}
Let $B=\bigcup_{\alpha=1}^N L_{\alpha}$ and
$\overline{B}=\bigcup_{\beta=1}^{\overline{N}} \overline{L}_{\beta}$ be two
infinity branches of two algebraic curves ${\cal C}$ and
$\overline{{\cal C}}$, respectively.
\begin{enumerate}
\item $B$ and $\overline{B}$ converge if there are two convergent leaves  $L_\alpha\subseteq B, \alpha=1,\ldots,N$ and
$\overline{L}_\beta\subseteq \overline{B}, \beta=1,\ldots,\overline{N}$.
\item $B$ and $\overline{B}$ diverge if any two leaves $L_\alpha\subseteq B, \alpha=1,\ldots,N$ and
$\overline{L}_\beta\subseteq \overline{B}, \beta=1,\ldots,\overline{N}$
diverge.
\end{enumerate}
\end{definition}

From Definition \ref{D-convergent-branches} we get that two infinity branches
$B$ and $\overline{B}$ do not diverge if there are two leaves,
$L\subseteq B$ and $\overline{L}\subseteq \overline{B}$, that do not
diverge. 
Furthermore, the next lemma states that, in this case, every leaf of
$B$ is non-divergent with some leaf of $\overline{B}$, and
reciprocally.

\para


\begin{lemma}\label{L-non-div-branches} Let $B=\bigcup_{\alpha=1}^N L_{\alpha}$ and
$\overline{B}=\bigcup_{\beta=1}^{\overline{N}} \overline{L}_{\beta}$  be two non-divergent infinity branches. Then, for each leaf $L_\alpha\subseteq B$ there exists a leaf
$\overline{L}_\beta\subseteq \overline{B}$ that does not diverge
with $L_\alpha$, and reciprocally.
\end{lemma}
\noindent\textbf{Proof:} Let $B$ and $\overline{B}$ be two
non-divergent branches. Let us prove that for any leaf
$L_\alpha\subseteq B$ there exist one or more leaves
$\overline{L}_\beta\subseteq \overline{B}$ non-divergent with
$L_\alpha$, and reciprocally. From the discussion above, we know that there exist
two leaves $\{(z,r_2(z),\ldots, r_n(z))\in {\Bbb C}^n:\,z\in {\Bbb
C},\,|z|>M\}\subset B$  and
$\{(z,\overline{r}_2(z),\ldots,\overline{r}_n(z))\in
{\Bbb C}^n:\,z\in {\Bbb C},\,|z|>\overline{M}\}\subset \overline{B}$
that do not diverge. Let
$$r_j(z)=z\varphi_j(z^{-1})=m_jz+u_{1,j}z^{1-\frac{N_{1,j}}{N}}+\cdots+u_{k,j}z^{1-\frac{N_{k,j}}{N}}+u_{k+1,j}z^{1-\frac{N_{k+1,j}}{N}}+\cdots,$$
$$\overline{r}_j(z)=z\overline{\varphi}_j(z^{-1})=\overline{m}_jz+\overline{u}_{1,j}z^{1-\frac{\overline{N}_{1,j}}{\overline{N}}}
+\cdots+\overline{u}_{k,j}z^{1-\frac{\overline{N}_{k,j}}{\overline{N}}}
+\overline{u}_{k+1,j}z^{1-\frac{\overline{N}_{k+1,j}}{\overline{N}}}+\cdots,$$
where $\overline{u}_{i,j} u_{i,j}\not=0$,
$N=\nu(B)=\lcm(N_2,\ldots,N_n)$,
$\overline{N}=\nu(\overline{B})=\lcm(\overline{N}_2,\ldots,\overline{N}_n)$,
$N_{k,j}< N \leq N_{k+1,j}$ and $\overline{N}_{k,j}< \overline{N}
\leq \overline{N}_{k+1,j}$ for some $k\in {\Bbb N}$ (note that $k$ may depend on $j$). Note also
that the expression above differs slightly from that of
(\ref{Eq-r_j}), since we are using $N$ and $\overline{N}$ as the
common denominators for the exponents of the series $r_j$ and
$\overline{r}_j$ respectively.

\para

From Lemma \ref{L-DistVertical}, we deduce that the terms with
positive exponent in $r_j$ and $\overline{r}_j$ are the same. Thus,
$\overline{m}_j=m_j$, $\overline{u}_{i,j}=u_{i,j}$, for
$i=1,\ldots,k$, $j=2,\ldots,n$, and
$$r_j(z)=m_jz+u_{1,j}z^{1-\frac{n_{1,j}}{n}}+\cdots+u_{k,j}z^{1-\frac{n_{k,j}}{n}}+u_{k+1,j}z^{1-\frac{N_{k+1,j}}{N}}+\cdots,$$
$$\overline{r}_j(z)=m_jz+u_{1,j}z^{1-\frac{n_{1,j}}{n}}+\cdots+u_{k,j}z^{1-\frac{n_{k,j}}{n}}+\overline{u}_{k+1,j}z^{1-\frac{\overline{N}_{k+1,j}}{\overline{N}}}
+\cdots, $$ where $\overline{u}_{i,j}, u_{i,j}\not=0$, $n,
n_{i,j}\in {\Bbb N}$ and $0<n_{1,j}<\cdots<n_{k,j}<n$. Observe that
we have simplified the non negative exponents such that
$\gcd(n,n_{1,j},\ldots,n_{k,j})=1$ ,for $j=2,\ldots,n$ . Hence, there
are $b,\overline{b}\in\mathbb{N}$ such that $N_{i,j}=bn_{i,j}$,
$N=bn$, $\overline{N}_{i,j}=\overline{b}n_{i,j}$, and
$\overline{N}=\overline{b}n$
 for $i=1,\ldots,k$ and $j=2,\ldots,n$.

\para

Under these conditions, we observe that the different leaves of $B$
and $\overline{B}$ are obtained by conjugation on $r_j(z)$ and
$\overline{r}_j(z)$, $j=2,\ldots,n$. That is, any two leaves
$L_{\alpha}\subseteq B,\, \alpha=1,\ldots,N$ and
$\overline{L}_{\beta}\subseteq \overline{B},\,
\beta=1,\ldots,\overline{N}$ will have the form
$L_{\alpha}=\{(z,r_{\alpha,2}(z),\ldots, r_{\alpha,n}(z))\in {\Bbb
C}^n:\,z\in {\Bbb C},\,|z|>M\}$ and
$\overline{L}_{\beta}=\{(z,\overline{r}_{\beta,2}(z),\ldots,\overline{r}_{\beta,n}(z))\in
{\Bbb C}^n:\,z\in {\Bbb C},\,|z|>\overline{M}\}$, where
$r_{\alpha,j}(z)=$  $$m_jz+u_{1,j}c_{\alpha}^{N_{1,j}}z^{1-\frac{N_{1,j}}{N}}+\cdots+u_{k,j}c_{\alpha}^{N_{k,j}}z^{1-\frac{N_{k,j}}{N}}+u_{k+1,j}c_{\alpha}^{N_{k+1,j}}z^{1-\frac{N_{k+1,j}}{N}}+\cdots,$$
and $\overline{r}_{\beta, j}(z)=$ $$\overline{m}_jz+\overline{u}_{1,j}d_{\beta}^{\overline{N}_{1,j}}z^{1-\frac{\overline{N}_{1,j}}{\overline{N}}}
+\cdots+\overline{u}_{k,j}d_{\beta}^{\overline{N}_{k,j}}z^{1-\frac{\overline{N}_{k,j}}{\overline{N}}}
+\overline{u}_{k+1,j}d_{\beta}^{\overline{N}_{k+1,j}}z^{1-\frac{\overline{N}_{k+1,j}}{\overline{N}}}+\cdots,$$
 $c_1,\ldots,c_{N}$ are the $N$ complex roots of $x^{N}=1$, and
$d_1,\ldots,d_{\overline{N}}$ are the $\overline{N}$ complex roots
of $x^{\overline{N}}=1$ (see equation (\ref{Eq-conjugates})).

\para

We simplify the exponents and, using that
$\overline{u}_{i,j}=u_{i,j},\, i=1,\ldots,k$, we get that:
$$r_{\alpha,j}(z)=m_jz+u_{1,j}c_{\alpha}^{N_{1,j}}z^{1-\frac{n_{1,j}}{n}}+\cdots+u_{k,j}c_{\alpha}^{N_{k,j}}z^{1-\frac{n_{k,j}}{n}}+u_{k+1,j}c_{\alpha}^{N_{k+1,j}}z^{1-\frac{N_{k+1,j}}{N}}+\cdots$$
$$\overline{r}_{\beta,j}(z)=m_jz+u_{1,j}d_{\beta}^{\overline{N}_{1,j}}z^{1-\frac{n_{1,j}}{n}}+\cdots+u_{k,j}d_{\beta}^{\overline{N}_{k,j}}z^{1-\frac{n_{k,j}}{n}}
+\overline{u}_{k+1,j}d_{\beta}^{\overline{N}_{k+1,j}}z^{1-\frac{\overline{N}_{k+1,j}}{\overline{N}}}+\cdots.$$

\para

Now, we   prove that  for any leaf $L_\alpha$ there exist one
or more leaves $\overline{L}_\beta$ non-divergent with $L_\alpha$.
For this purpose, we just need to show that, given any value of
$\alpha=1,\ldots,N$, there exist one or more values of
$\beta=1,\ldots,\overline{N}$ such that
$c_{\alpha}^{N_{i,j}}=d_{\beta}^{\overline{N}_{i,j}}, i=1,\ldots,k,
j=2,\ldots,n$.

\para

Indeed, since the coefficients $c_\alpha,\,\alpha=1,\ldots,N$ are
the $N$ complex roots of $x^{N}=1$, we have that
$c_\alpha=e^{\frac{2(\alpha-1)\pi I}{N}}$, where $I$ is the
imaginary unit. Taking into account that $N=bn$, we deduce that
$c_\alpha^b=e^{\frac{2(\alpha-1)\pi I}{n}}$ for each
$\alpha=1,\ldots,N$ and $c_\alpha^b=c_{\alpha+(m-1)n}^b$ for each
$\alpha=1,\ldots,n$ and $m=1,\ldots,b$. That is,
$c_\alpha^b,\,\alpha=1,\ldots,n$ are the $n$ complex roots of
$x^{n}=1$. Reasoning similarly, we have that
$d_\beta^{\overline{b}}=e^{\frac{2(\beta-1)\pi I}{n}}$ for each
$\beta=1,\ldots,\overline{N}$ and
$d_\beta^{\overline{b}}=d_{\beta+(m-1)n}^{\overline{b}}$ for each
$\beta=1,\ldots,n$ and $m=1,\ldots,\overline{b}$. That is,
$d_\beta^{\overline{b}},\,\beta=1,\ldots,n$ are the $n$ complex
roots of $x^{n}=1$. Hence, for each $\alpha=1,\ldots,N$ there are
one or more $\beta=1,\ldots,\overline{N}$ such that
$c_\alpha^b=d_\beta^{\overline{b}}$, and reciprocally. Finally, the
result follows  taking into account that
$c_\alpha^{N_{i,j}}=\left(c_\alpha^b\right)^{n_{i,j}}=\left(d_\beta^{\overline{b}}\right)^{n_{i,j}}=d_\beta^{\overline{N}_{i,j}}$.
\hfill $\Box$

\begin{remark}
Let $B$  and $\overline{B}$ be two infinity branches associated with
two infinity points $P=(1:m_2:\cdots:m_n)$ and
$\overline{P}=(1:\overline{m}_2:\cdots:\overline{m}_n)$,
respectively. From the proof of Lemma \ref{L-non-div-branches}, if
$B$  and $\overline{B}$ do not diverge, then $m_j=\overline{m}_j$
for every $j=2,\ldots,n$ which implies that two non-divergent infinity
branches are  associated with the same infinity point (see
Remark 4.5 in \cite{paper1}).
\end{remark}

\para

For the sake of simplicity, and taking into account that an infinity branch $B$  is uniquely
determined  from one leaf, up to conjugation (see
Remark \ref{R-conjugation}),  we identify an infinity branch by
just one of its leaves. Hence, in the following
 $$B=\{(z,r_2(z),\ldots, r_{n}(z))\in {\Bbb C}^n:\,z\in {\Bbb
C},\,|z|>M\},\qquad M\in {\Bbb R}^+$$
will stand for the infinity branch whose leaves are obtained by
conjugation on
$$r_{j}(z)=
m_jz+a_{1,j} z^{1-N_{1,j}/N_j}+a_{2,j} z^{1-N_{2,j}/N_j}+a_{3,j}
z^{1-N_{3,j}/N_j} + \cdots,$$ $a_{i,j}\not=0,\, \forall i\in {\Bbb
N},\,i\geq 1,$ $N_j, N_{i,j}\in {\Bbb N}$, and
$0<N_{1,j}<N_{2,j}<\cdots$ for $j=2,\ldots, n$. Observe that   the
results stated above hold for any leaf of $B$.

\para

 Finally, we remark that there exists well known algorithms that allow to compute the series
 $\varphi_{j}(t)\in {{\Bbb C}\ll t\gg},\,j=2,\ldots,n$, and then the branch $B=\{(z,r_2(z),\ldots, r_{n}(z))\in {\Bbb C}^n:\,z\in {\Bbb
C},\,|z|>M\}$  (see e.g. \cite{Alon92}). In addition, in
\cite{paper3}, a procedure  for computing the branches for $n=3$ is
presented. This method is based on projections over the plane, and
it can be generalized for a given curve in the $n$-dimensional space
by successively eliminating variables  and reducing the problem to
the computation of infinity branches for plane curves (a method for
successively eliminating the variables, by means of univariate
resultants, is presented in \cite{PDS-Fibra}). For the plane case
($n=2$) methods are well known (see e.g. \cite{BlascoPerezII},
\cite{paper1}).

\para

In the following example, we compute the infinity branches for a given algebraic curve  in the $4$-dimensional
space  implicitly defined by the polynomials $f_i(x_1,x_2,x_3,x_4)\in {\Bbb R}[x_1,x_2,x_3,x_4],\,i=1,2,3$.

\para

\begin{example}\label{E-infbranches-impli}
Let $\cal C$ be the irreducible curve defined over $\Bbb C$ by the
polynomials $$f_1(x_1,x_2,x_3,x_4)=x_1-x_2^2+2x_3 ,\quad
f_2(x_1,x_2,x_3,x_4)=x_1+x_2-x_4^2 ,\,\,\mbox{ and}\,\,\,\,\,$$$$
f_3(x_1,x_2,x_3,x_4)=2x_2-x_3^2+x_4.$$ The projection along the
$x_4$-axis, ${\cal C}^p$,  is defined by the polynomials
$$f_1^p(x_1,x_2,x_3)=x_1-x_2^2+2x_3,\quad \mbox{and}\quad f_2^p(x_1,x_2,x_3)=x_1+x_2-4x_2^2+4x_2x_3^2-x_3^4$$ (these polynomials can be
obtained by computing  univariate resultants).
\para
By applying the method described in \cite{paper3}, we compute
the infinity branches of ${\cal C}_p$. We obtain the branch
$B_1^p=\{(z,r_{1,2}(z),r_{1,3}(z))\in {\Bbb C}^3:\,z\in {\Bbb
C},\,|z|>M_1^p\}$, where
$$r_{1,2}(z)=z^{1/2}+\sqrt{3}z^{-1/4}+\frac{\sqrt{3}z^{-3/4}}{12}-\frac{z^{-1}}{2}-\frac{7\sqrt{3}z^{-5/4}}{288}+\cdots$$
$$r_{1,3}(z)=\sqrt{3}z^{1/4}+\frac{\sqrt{3}z^{-1/4}}{12}+z^{-1/2}-\frac{7\sqrt{3}z^{-3/4}}{288}+\frac{z^{-1}}{4}+\cdots,$$
and the branch $B_2^p=\{(z,r_{2,2}(z),r_{2,3}(z))\in {\Bbb C}^3:\,z\in
{\Bbb C},\,|z|>M_2^p\}$, where
$$r_{2,2}(z)=z^{1/2}+z^{-1/4}-\frac{z^{-3/4}}{4}+\frac{z^{-1}}{2}-\frac{z^{-5/4}}{32}+\cdots,$$
$$r_{2,3}(z)=z^{1/4}-\frac{z^{-1/4}}{4}+z^{-1/2}+\frac{z^{-3/4}}{32}-\frac{z^{-1}}{4}+\cdots.$$
Note that both branches are associated to the infinity point
$P_1=(1:0:0:0)$. Moreover, $\nu(B_1^p)=\nu(B_2^p)=4$, and thus each branch has 4 (conjugated) leaves. That is, $B_1^p=\bigcup_{\alpha=1}^4 L_{1,\alpha}$, where $L_{1,\alpha}$ are obtained by conjugation in the above series $r_{1,2}$ and $r_{1,3}$ (similarly  for $B_2^p$).

\para

\noindent Once we have the infinity branches of the projected curve
${\cal C}^p$, we compute the infinity branches of the curve ${\cal
C}$. We use the lift function $h(x_1,x_2,x_3)=-2x_2+x_3^2$ to get
the fourth component of these branches (we apply the results in
\cite{Bajaj} to compute $h$). Thus, the infinity branches of the
curve $\cal C$ are $B_1=\{(z,r_{1,2}(z),r_{1,3}(z),r_{1,4}(z))\in {\Bbb
C}^4:\,z\in {\Bbb C},\,|z|>M_1\}$, where
$$r_{1,4}(z)=h(z,r_{1,2}(z),r_{1,3}(z))=z^{1/2}+\frac{1}{2}-\frac{z^{-1/2}}{8}+\frac{\sqrt{3}z^{-3/4}}{2}+\cdots$$ and
$B_2=\{(z,r_{2,2}(z),r_{2,3}(z),r_{2,4}(z))\in {\Bbb C}^4:\,z\in {\Bbb
C},\,|z|>M_2\}$, where
$$r_{2,4}(z)=h(z,r_{2,2}(z),r_{2,3}(z))=-z^{1/2}-\frac{1}{2}+\frac{z^{-1/2}}{8}-\frac{z^{-3/4}}{2}+\cdots.$$
In Figure \ref{F-ejemplo-ramas-implicitas}, we plot the curve ${\cal
C}^p$ and some points of the infinity branches $B_1^p$ and $B_2^p$.

\begin{figure}[h]
$$
\begin{array}{cc}
\psfig{figure=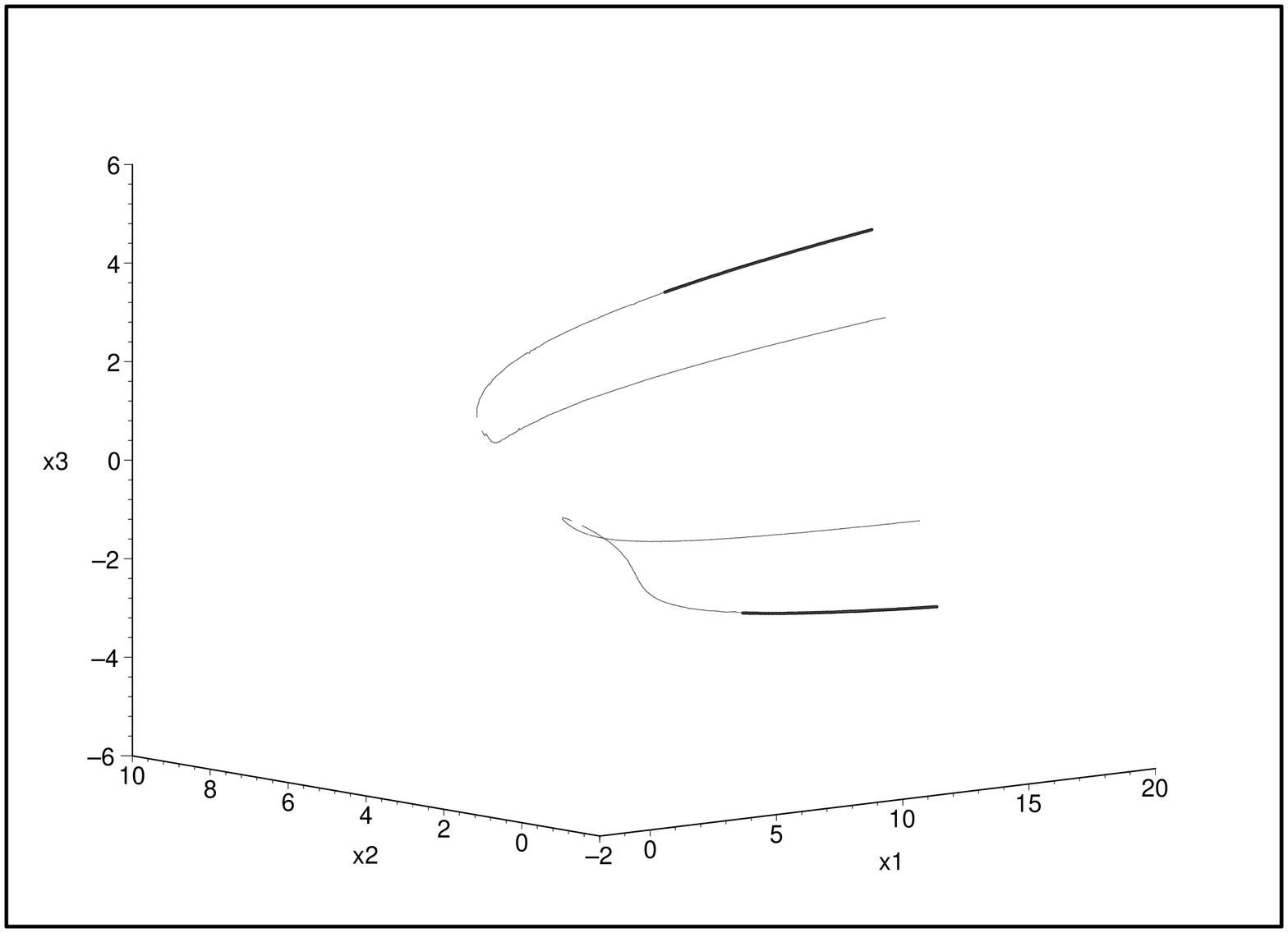,width=5cm,height=5cm,angle=270} &
\psfig{figure=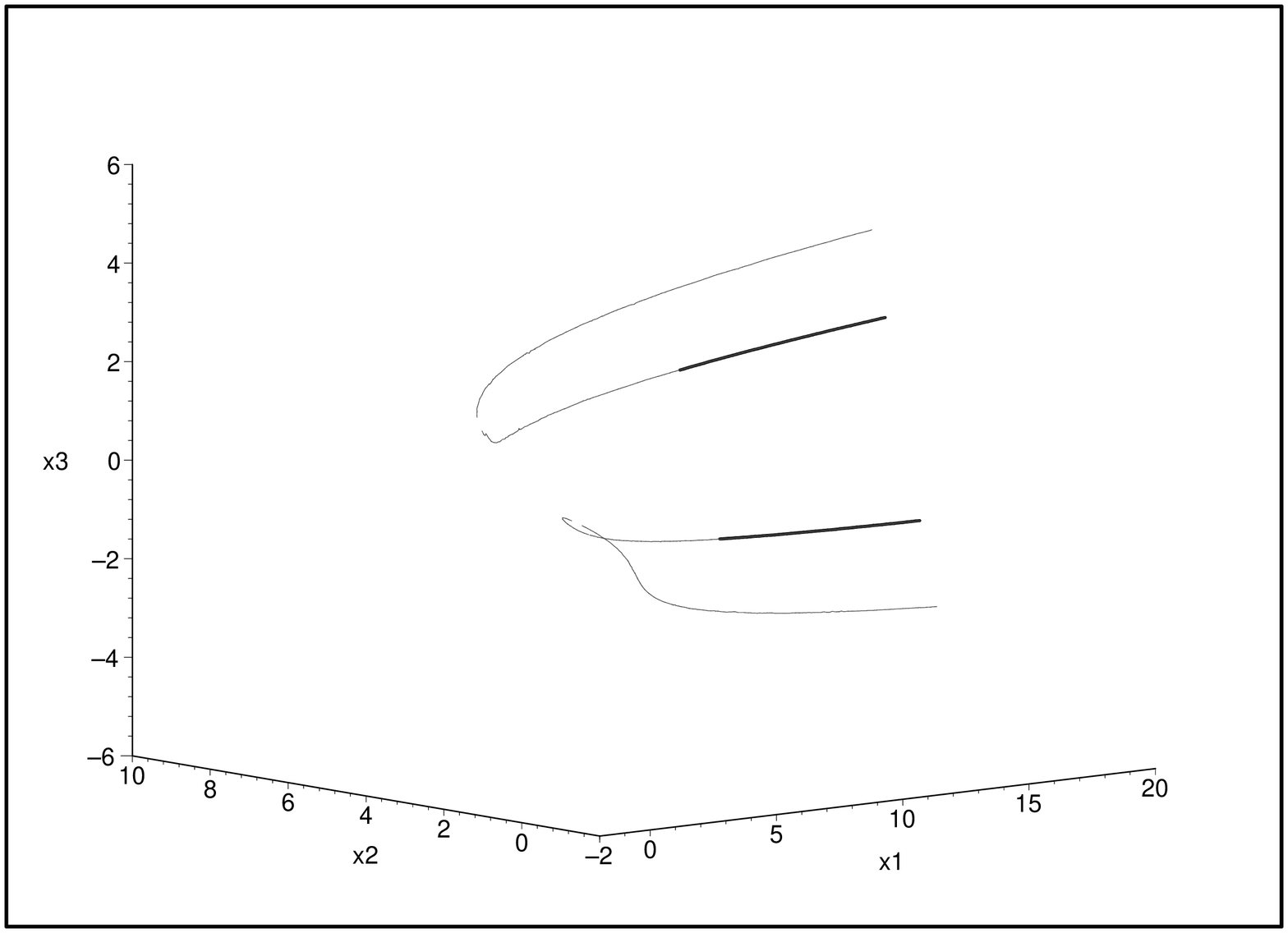,width=5cm,height=5cm,angle=270}
\end{array}
$$ \caption{Curve ${\cal
C}^p$ and infinity branches $B_1^p$ (left) and $B_2^p$
(right).}\label{F-ejemplo-ramas-implicitas} \vspace*{-0.05cm}
\end{figure}
\end{example}

\subsection{Parametrically defined space curves} \label{Sub-Param}

In Subsection \ref{Sub-Impli}, we have assumed that the given real algebraic  curve  in the $n$-dimensional space is defined
implicitly by some polynomials.  In this section, we show how to deal with rational curves defined parametrically.

\para

Note that the definitions introduced above are independent on whether the curve is defined parametrically or implicitly. However, the method to compute the infinity branches has to be different (of course, one may implicitize and reason as in  Subsection \ref{Sub-Impli}, but we are interested in computing the infinity branches from the given parametrization without implicitizing).

\para

Thus, in this subsection, we present a method
to compute infinity branches  of a rational  curve  in the $n$-dimensional space  from
their parametric representation (without implicitizing).
Similarly as above,   we  work
over   $\Bbb C$, but we assume that the curve has infinitely many points in the affine plane over $\Bbb R$ and then, the curve has a real parametrization. The method presented generalize  the results in \cite{paper3}.

\para

Under these conditions, in the following, we consider a real space
curve ${\cal C}$ in the $n$-dimensional space $\mathbb{C}^n$,
defined by the parametrization
\[{\cal P}(s)=(p_1(s),\ldots,p_n(s))\in {\Bbb R}(s)^n\setminus{\Bbb R}^n,\quad p_i(s)=p_{i1}(s)/p(s),\,i=1,\ldots,n.\]
We assume that   we have  prepared the input curve ${\cal C}$, by
means of a suitable linear change of coordinates (if necessary) such
that $(0:m_2:\ldots:m_n:0)$ ($m_j\not=0$ for some $j=2,\ldots,n$) is
not an infinity point  (see Remark \ref{R-infinitypoint}).  Note that, hence, $\deg(p_1)\geq 1$.

\para

Now, let ${\cal C}^*$ denote  the projective curve associated to
${\cal C}$. We have that a parametrization of ${\cal C}^*$ is given
by ${\cal P}^*(s)=(p_{11}(s):\cdots:p_{n1}(s):p(s))$ or,
equivalently,
$${\cal P}^*(s)=\left(1:\frac{p_{21}(s)}{p_{11}(s)}:\cdots:\frac{p_{n1}(s)}{p_{11}(s)}:\frac{p(s)}{p_{11}(s)}\right).$$

\para

Under these conditions, we show how  to compute the   infinity
branches of $\cal C$. That is, the sets
$B=\{(z:r_{2}(z):\ldots:r_{n}(z)):\,z\in {\Bbb C},\,|z|>M\},$ where
$r_{j}(z)=z\varphi_{j}(z^{-1})\in{{\Bbb C}\ll z\gg},\,
j=2,\ldots,n$. We recall  that these
series must verify $F_i(1:\varphi_{2}(t):\ldots:\varphi_{n}(t):t)=0$
around $t=0$, where $F_i,\,i=1,\ldots,s$ are the polynomials
defining implicitly ${\cal C}^*$ (see Subsection \ref{Sub-Impli}).  Observe that in this subsection,
we are given the parametrization ${\cal P}^*$ of ${\cal C}^*$ and
then, $F_i({\cal
P}^*(s))=F_i \left(1:\frac{p_{21}(s)}{p_{11}(s)}:\cdots:\frac{p_{n1}(s)}{p_{11}(s)}:\frac{p(s)}{p_{11}(s)}\right)=0$.
Thus, intuitively speaking, in order to compute the  infinity
branches of $\cal C$, and in particular the series
$\varphi_{j},\,j=2,\ldots,n$, one needs to ``reparametrize'' the
parametrization ${\cal
P}^*(s)=\left(1:\frac{p_{21}(s)}{p_{11}(s)}:\ldots:\frac{p_{n1}(s)}{p_{11}(s)}:\frac{p(s)}{p_{11}(s)}\right)$
in the form $(1:\varphi_{2}(t):\ldots:\varphi_{n}(t):t)$ around
$t=0$. For this purpose, the idea is to look for a value of the
parameter $s$, say $\ell(t)\in {\Bbb C}\ll
t\gg$, such that ${\cal P}^*(\ell(t))=(1:\varphi_{2}(t):\ldots:\varphi_{n}(t):t)$ around $t=0$.\\

Hence, from the above reasoning, we deduce that first, we have to
consider the equation $p(s)/p_{11}(s)=t$ (or equivalently,
$p(s)-tp_{11}(s)=0$), and we have to solve it in the variable $s$
around $t=0$ (note that $\deg(p_1)\geq 1$). From
Puiseux's Theorem, there exist solutions
$\ell_1(t),\ell_2(t),\ldots,\ell_k(t)\in {\Bbb C}\ll t\gg$,
where $k=\deg(p_1)$, such that,
$p(\ell_i(t))-tp_{11}(\ell_i(t))=0,\,i=1,\ldots,k,$ in a neighborhood of $t=0$.\\

Thus,   for each $i=1,\ldots,k$, there exists  $M_i\in {\Bbb R}^+$ such that
the points $(1:\varphi_{i,2}(t):\ldots:\varphi_{i,n}(t):t)$ or  equivalently, the
points $(t^{-1}:t^{-1}\varphi_{i,2}(t):\ldots:t^{-1}\varphi_{i,n}(t):1)$, where
\begin{equation}\label{Eq-psi}\varphi_{i,j}(t)=\frac{p_{j,1}(\ell_i(t))}{p_{11}(\ell_i(t))},\quad
j=2,\ldots,n,\end{equation} are in ${\cal C}^*$ for $|t|<M_i$.
Observe  that $\varphi_{i,j}(t),\, j=2,\ldots,n,$ are
Puiseux series, since $p_{j,1}(\ell_i(t)),\, j=2,\ldots,n$, and
$p_{11}(\ell_i(t))$ can be written as Puiseux series (around $t=0$) and ${\Bbb C}\ll t\gg$ is a
field.\\

Finally, we set $z=t^{-1}$. Then, we have that the points
$(z:r_{i,2}(z):\ldots:r_{i,n}(z))$, where $r_{i,j}(z)=z\varphi_{i,j}(z^{-1}),
j=2,\ldots,n$, are in ${\cal C}$ for $|z|>M_i^{-1}$. Hence, the infinity
branches of $\cal C$ are the sets
$$B_i=\{(z:r_{i,2}(z):\ldots:r_{i,n}(z))\in {\Bbb C}^n:\,z\in {\Bbb
C},\,|z|>M_i^{-1}\},\quad i=1,\ldots,k.$$

\para

\begin{remark}\label{R-calculo-r-directo} We observe that:
\begin{enumerate}
\item The series $\ell_i(t)$ satisfies that
$p(\ell_i(t))/p_{11}(\ell_i(t))=t$, for $i=1,\ldots,k$. Then, from equality
(\ref{Eq-psi}), we have that for $j=2,\ldots,n$
$$\varphi_{i,j}(t)=\frac{p_{j,1}(\ell_i(t))}{p(\ell_i(t))}t=p_j(\ell_i(t))t,\,\,\mbox{and}\,\, r_{i,j}(z)=z\varphi_{i,j}(z^{-1})=p_j(\ell_i(z^{-1})).$$
\item In order to compute $r_{i,j}(z)$, we first write    $p_j(\ell_i(t))$  as Puiseux series  around $t=0$, and then we set $t=z^{-1}$.
\item  When we compute the series ${\ell}_i$, we cannot handle its infinite terms so it must be truncated,
which may distort the computation of the series $r_{i,j}$. The
number of affected terms in $r_{i,j}$ depends on the number of terms
computed in ${\ell}_i$. That is, as more terms we compute in
${\ell}_i$, as more accurate the computation of $r_{i,j}$ is. More
details on this question are analyzed in Proposition 5.4 in
\cite{paper3}.\end{enumerate}
\end{remark}

\para

In the following example, we show the above procedure and we compute the infinity branches for a given  curve defined by a parametrization $\cP(s)\in {\Bbb R}(s)^4$.

\para

\begin{example}\label{E-infbranches-para}
Let ${\cal C}$ be the  curve defined by the
parametrization
$${\cal P}(s)=\left(p_1(s), p_2(s), p_3(s), p_4(s)\right)=\left(\frac{p_{11}(s)}{p(s)},\,\frac{p_{21}(s)}{p(s)},\,\frac{p_{31}(s)}{p(s)},\,\frac{p_{41}(s)}{p(s)}\right)=$$
$$=\left(\frac{-1+2s^3-s}{s},\,\frac{s+1}{s},\,\frac{-1}{s},\,\frac{s^2+3s-5}{s}\right)\in {\Bbb R}(s)^4.$$
 We compute the solutions of the equation $p(s)-tp_{11}(s)=0$
around $t=0$. We get the
Puiseux series
$$\ell_1(t)=-t+t^2-t^3-t^4+7t^5+\cdots$$
$$\ell_2(t)=\frac{1}{2}\sqrt{2}t^{-1/2}+\frac{1}{4}\sqrt{2}{t}^{1/2}+\frac{1}{2} t-\frac{1}{16}\sqrt{2} t^{3/2}-\frac{1}{2} t^2-\frac{11}{32} \sqrt{2}t^{5/2}+\frac{1}{2} t^3+\frac{235}{256}  \sqrt{2} t^{7/2}+\cdots$$
(note that $\ell_2(t)$ represents a conjugation class composed by two  conjugated
series).\\

\noindent
 Now, we determine the series $r_{i,j}(z),\,i=1,2,\,j=2,3,4$. We get
 $$\begin{array}{l}
 r_{1,2}(z)=p_2(\ell_1(z^{-1}))=-z+2 z^{-2}-4z^{-3}-13 z^{-4}-11 z^{-5}    +\cdots\\
 \\
 r_{1,3}(z)=p_3(\ell_1(z^{-1}))=z+1-2z^{-2} +4z^{-3} +13z^{-4} +11z^{-5}     +\cdots\\
 \\
 r_{1,4}(z)=p_4(\ell_1(z^{-1}))=5z+8-z^{-1} -9z^{-2} +19z^{-3} +64z^{-4} +62z^{-5}    +\cdots,
 \end{array}$$
 and
 $$\begin{array}{l}
 r_{2,2}(z)=p_2(\ell_2(z^{-1}))=1+\sqrt{2}z^{-1/2} -\frac{1}{2}\sqrt{2}z^{-3/2}-z^{-2}+\frac{3}{8}\sqrt{2}z^{-5/2}+2z^{-3}+\cdots\\
 \\
 r_{2,3}(z)=p_3(\ell_2(z^{-1}))=-\sqrt{2}z^{-1/2}+\frac{1}{2}\sqrt{2}z^{-3/2}+z^{-2}-\frac{3}{8}\sqrt{2}z^{-5/2}-2z^{-3}+\cdots\\
 \\
 r_{2,4}(z)=p_4(\ell_2(z^{-1}))=\frac{1}{2}\sqrt{2}z^{1/2} +3-\frac{19}{4}\sqrt{2}z^{-1/2}+\frac{1}{2}z^{-1}+\frac{39}{16}\sqrt{2}z^{-3/2}+\\
 \qquad \qquad \qquad \qquad  \qquad \frac{9}{2} z^{-2}-\frac{71}{32}\sqrt{2}z^{-5/2}-\frac{19}{2} z^{-3}+\cdots\\
 \end{array}$$

 \noindent
Therefore, the curve has two infinity branches given by
$$B_1=\{(z,r_{1,2}(z),r_{1,3}(z),r_{1,4}(z))\in {\Bbb C}^4:\,z\in {\Bbb
C},\,|z|>M_1\}$$
and
$$B_2=\{(z,r_{2,2}(z),r_{2,3}(z),r_{2,4}(z))\in {\Bbb C}^4:\,z\in {\Bbb
C},\,|z|>M_2\}$$
for some $M_1, M_2\in {\Bbb R}^+$. Note that $B_1$ is associated   to the infinity point $(1:-1:1:5:0)$, and  $B_2$ is associated   to the infinity point $(1:0:0:0:0)$. In addition, we observe that $\nu(B_1)=1$ and $\nu(B_2)=2$,  and thus $B_1$ has one leaf, and  $B_2$ has two (conjugated) leaves.
\end{example}

\section{Asymptotic behavior and Hausdorff distance}

In this section, we consider algebraic  curves in the $n$-dimensional space
defined by a finite set of real polynomials or by  a rational parametrization. Depending on whether they are defined parametrically or implicitly
one proceeds as in Subsection \ref{Sub-Impli} or as in Subsection
\ref{Sub-Param} to compute their infinity branches.

\para

We remind that the input curves are prepared such that
$(0:m_2:\ldots: m_n:0)$ ($m_j\not=0$ for some $j=2,\ldots,n$) is not
an infinity point of their corresponding projective curves (see
Remark \ref{R-infinitypoint}).

\para

The main result of the section states that the Hausdorff distance
between two algebraic curves is finite if and only if their
asymptotic behaviors are similar (we say that two algebraic curves
have similar asymptotic behaviors if their infinity branches are
pair-wise non-divergent; see Definition \ref{D-sim-behavior}).

\para

The computation of the Hausdorff distance plays an important role in
the frame of practical applications in computer aided geometric design such as approximate
parametrization problems (see Section 1). In particular, estimating
the Hausdorff distance between two curves is specially interesting
since it is an appropriate tool for measuring the closeness between
them. Many authors have addressed some problems in this frame  (see
e.g. \cite{Bai}, \cite{Chen}, \cite{Kim}, \cite{Patri},  \cite{RSS},
etc).\\

\noindent To start with, we first introduce the following
definition.

\para

\begin{definition}\label{D-sim-behavior}
We say that two algebraic curves, ${\cal C}$ and ${\overline{{\cal
C}}}$, have a {\sf similar asymptotic behavior} if, for every
infinity branch $B\subseteq{\cal C}$ there exist an infinity branch
$\overline{B}\subseteq{\overline{{\cal C}}}$ non-divergent with $B$,
and reciprocally.
\end{definition}

\para\para

  Now, we introduce the notion of  {\em Hausdorff distance}. For this purpose, we recall that,
  given an algebraic space curve ${\cal C}$ over $\Bbb C$ and a point $p\in {\Bbb
C}^n$,   {\em the distance from $p$ to ${\cal C}$} is defined as
$d(p,{\cal C})=\min\{d(p,q):\,q\in{\cal C}\}.$



\para

\begin{definition}\label{D-Hausdorff} Given a metric space $(E,d)$ and two
subsets $A, B\subset E\setminus \{\emptyset\}$, the {\em Hausdorff distance} between them is defined
as:
$$d_H(A,B)=\max\{\sup_{x\in A}\inf_{y\in B}d(x,y),\sup_{y\in B}\inf_{x\in A}d(x,y)\}.$$
If $E=\mathbb{C}^{n}$ and $d$ is the unitary  distance, the Hausdorff distance between two curves ${\cal C}$ and
${\overline{{\cal C}}}$ can be expressed as:
$$d_H({\cal C},{\overline{{\cal C}}})=\max\{\sup_{p\in {\cal C}}d(p,{\overline{{\cal C}}}),\sup_{\overline{p}\in {\overline{{\cal C}}}}d(\overline{p},{\cal C})\}.$$
\end{definition}

\para

In order to prove the main theorem (see Theorem \ref{P-CharacterizationHausdorff}), we first need to prove some
technical lemmas. The first one  (Lemma \ref{L-BolaMaxima}) states
that  any point of the curve with sufficiently large coordinates
belongs to some infinity branch (see also Lemma 3.6 and
Remark 3.7 in \cite{paper1}).

 \para

\begin{lemma}\label{L-BolaMaxima} Let $\cal C$ be an algebraic space curve.
There exists $K\in\mathbb{R}^+$ such that every $p=(a_1,\ldots,
a_n)\in\mathcal{C}$ with $|a_i|>K$ (for some $i\in \{1,\ldots, n\}$)
belongs to some infinity branch of ${\cal C}$.
\end{lemma}

\noindent\textbf{Proof:}  First, let us prove that there exists
$K^1\in\mathbb{R}^+$ such that every point
$p=(a_1,\ldots, a_n)\in\mathcal{C}$ with $|a_1|>K^1$ belongs to some
infinity branch.

\para

Let us assume that this is not true and let us consider a sequence
$\{K_{\kappa}\}_{{\kappa}\in {\Bbb N}}\in \mathbb{R}^+$ such
that $\lim_{{\kappa}\rightarrow\infty}K_{\kappa}=\infty$. Then, for
every ${\kappa}\in {\Bbb N}$ there   exists a point
$p_{\kappa}=(a_{1,{\kappa}},\ldots, a_{n,{\kappa}})\in\mathcal{C}$
such that $|a_{1,{\kappa}}|>K_{\kappa}$, and $p_{\kappa}$ does not
belong to any infinity branch of $\cal C$. The corresponding
projective point is $P_{\kappa}=(a_{1,{\kappa}}:\ldots
:a_{n,{\kappa}}:1)$, and it holds that $F(P_{\kappa})=f(p_{\kappa})=0$.
Thus, we have a sequence $\{P_{\kappa}\}_{{\kappa}\in {\Bbb N}}$ of
points in the projective curve ${\cal C}^*$ such that
$\lim_{{\kappa}\rightarrow\infty}|a_{1,{\kappa}}|=\infty$. Note that these projective points can be expressed as
$$P_{\kappa}=(1:a_{2,{\kappa}}/a_{1,{\kappa}}:\ldots
:a_{n,{\kappa}}/a_{1,{\kappa}}:1/a_{1,{\kappa}}).$$

Under these conditions, we extract
a subsequence $\{P_{\kappa_l}\}_{l\in {\Bbb N}}$ for the sequences
$\{a_{i,{\kappa_l}}/a_{1,{\kappa_l}}\}_{l\in {\Bbb N}}$,
$i=2,\ldots,n$ to be monotone. In order to simplify the notation, we also denote  it as $\{P_{\kappa}\}_{{\kappa}\in {\Bbb N}}$.
Now, we distinguish two different cases:
\begin{enumerate}
\item Let us assume that all these monotone sequences are bounded. Then,
$\lim_{{\kappa}\rightarrow\infty}a_{i,{\kappa}}/a_{1,{\kappa}}=m_i\in
{\Bbb C},\,i=2,\ldots,n$ and
$\lim_{{\kappa}\rightarrow\infty}1/a_{1,{\kappa}}=0$. Furthermore,
since $F(P_{\kappa})=0$ for every ${\kappa}\in {\Bbb N}$, we get that
$\lim_{{\kappa}\rightarrow\infty}F(P_{\kappa})=F(\lim_{{\kappa}\rightarrow\infty}P_{\kappa})=F(1:m_2:\cdots:m_n:0)=0$.
We conclude that the sequence $\{P_{\kappa}\}_{{\kappa}\in {\Bbb
N}}$ converges to the infinity point $P=(1:m_2:\cdots:m_n:0)$ as
${\kappa}$ tends to infinity; that is, there exists $M\in {\Bbb
R}^+$ such that $\|P_{\kappa}-P\|\leq \epsilon$, for ${\kappa}\geq
M$. Thus, we deduce that the points $\{P_{\kappa}\}_{{\kappa}\in
{\Bbb N},\,{\kappa}\geq M}$ can be obtained by a place centered at
$P$. Hence, the points $\{p_{\kappa}\}_{{\kappa}\in {\Bbb
N},\,{\kappa}\geq M}$ belong to some infinity branch of ${\cal C}$,
which contradicts the hypothesis.

\item If not all the sequences are bounded, then there is some
$i=2,\ldots,n$ such that
$\lim_{l\rightarrow\infty}a_{i,{\kappa}}/a_{1,{\kappa}}=\pm\infty$.
We assume without lost of generality that
$\lim_{l\rightarrow\infty}a_{2,{\kappa}}/a_{1,{\kappa}}=\pm\infty$.
 Then, we write
$$P_{\kappa}=(a_{1,{\kappa}}/a_{2,{\kappa}}:1:a_{3,{\kappa}}/a_{2,{\kappa}}:\ldots
:a_{n,{\kappa}}/a_{2,{\kappa}}:1/a_{2,{\kappa}}),$$ and we extract a
subsequence $\{P_{\kappa_l}\}_{l\in {\Bbb N}}$ for the sequences
$\{a_{i,{\kappa_l}}/a_{2,{\kappa_l}}\}_{l\in {\Bbb N}}$,
$i=3,\ldots,n$ to be monotone. For the sake of simplicity, we denote it by
$\{P_{\kappa}\}_{{\kappa}\in {\Bbb N}}$.

\para

At this point, we consider two different situations:
\begin{itemize}
\item If all these monotone sequences are bounded,  we get that
$$\lim_{{\kappa}\rightarrow\infty}a_{i,{\kappa}}/a_{1,{\kappa}}=m_i\in
{\Bbb C},\,i=3,\ldots,n.$$ Furthermore,
$\lim_{{\kappa}\rightarrow\infty}a_{1,{\kappa}}/a_{2,{\kappa}}=\lim_{{\kappa}\rightarrow\infty}1/a_{2,{\kappa}}=0$
 and thus,  reasoning as above, we deduce that the sequence
$\{P_{\kappa}\}_{{\kappa}\in {\Bbb N}}$ converges to an infinity
point $P=(0:1:m_3:\cdots:m_n:0)$.

\item If some of the sequences
$\{a_{i,{\kappa_l}}/a_{2,{\kappa_l}}\}_{l\in {\Bbb N}}$,
$i=3,\ldots,n$ are not bounded, we can assume w.l.o.g. that
$\lim_{l\rightarrow\infty}a_{3,{\kappa}}/a_{2,{\kappa}}=\pm\infty$
and we reason as above. Finally, we obtain a subsequence that converges
to an infinity point of the form $(0:m_2:m_3:\cdots:m_n:0)$.
\end{itemize}

In both cases, we find a contradiction, since we have prepared the input
curve such that it does not have infinity points of the form
$(0:m_2:m_3:\cdots:m_n:0)$.
\end{enumerate}

From the above discussion, the initial assumption leads us to a
contradiction. Therefore,  there exists
$K^1\in\mathbb{R}^+$ such that every point of the curve
$p=(a_1,\ldots, a_n)$ with $|a_1|>K^1$ belongs to some infinity
branch. Reasoning similarly, we deduce that for each
$i=2,\ldots,n$, there exists $K^i\in\mathbb{R}^+$ such that every
point of the curve $p=(a_1,\ldots, a_n)$ with $|a_i|>K^i$ belongs to
some infinity branch. Finally, the result follows by taking
$K=\min\{K^1,\ldots,K^n\}$.
 \hfill $\Box$

\para
\para

The following technical lemma states that, given two divergent
branches $B$ and $\overline{B}$, we can find points in $B$ as far as
we want from any point in $\overline{B}$ (and
reciprocally).

\para

\begin{lemma}\label{L-divergent-branches} Let $B=\{(z,r_2(z),\ldots, r_{n}(z))\in {\Bbb C}^n:\,z\in {\Bbb
C},\,|z|>M\}$ and $\overline{B}=\{(z,\overline{r}_2(z),\ldots, \overline{r}_{n}(z))\in {\Bbb C}^n:\,z\in {\Bbb
C},\,|z|>\overline{M}\}$ be two divergent infinity
branches. For each $K>0$, there exists $\delta>0$ such that if
$|x|>\delta$ then $d((x,r_2(x),\ldots, r_{n}(x)),(y,\overline{r}_2(y),\ldots, \overline{r}_{n}(y)))>K$ for any
point $(y,\overline{r}_2(y),\ldots, \overline{r}_{n}(y))\in \overline{B}$.
\end{lemma}
\noindent\textbf{Proof:} We assume w.l.o.g.  that $B$ is associated
to the infinity point $(1:0:\ldots :0)$ (otherwise we can apply a
linear change of coordinates). Note that since all the norms in
$\mathbb{C}^n$ are equivalent, there exists some $\lambda>0$ such
that
$$d((x,r_2(x),\ldots, r_{n}(x)),(y,\overline{r}_2(y),\ldots, \overline{r}_{n}(y)))>$$$$\lambda(|x-y|+|r_2(x)-\overline{r}_2(y)|+\cdots + |r_n(x)-\overline{r}_{n}(y)|).$$
Thus, we only need to prove that, for each $K>0$ there exists $\delta>0$ such that if
$|x|>\delta$ then $$\phi(x,y):=|x-y|+|r_2(x)-\overline{r}_2(y)|+\cdots + |r_n(x)-\overline{r}_{n}(y)|>K.$$
\noindent First of all, if $|x-y|>K$ the result follows, so we
assume that $|x-y|\leq K$. Hence, $|y|>|x|-K$ since $|x-y|>|x|-|y|$.\\
\noindent On the other hand, note that
$$|r_i(x)-\overline{r}_i(y)|=|\overline{r}_i(y)-r_i(x)|>|\overline{r}_i(y)-r_i(y)+r_i(y)-r_i(x)|>$$
\begin{equation}\label{Eq-ultimo-lema}>|\overline{r}_i(y)-r_i(y)|-|r_i(y)-r_i(x)|,\quad i=2,\ldots,n
\end{equation}
From the proof of Theorem 4.11 in \cite{paper1}, we get that
$r_i(z)$ is derivable for $|z|>M$ and $\limit_{z\rightarrow
\infty}r_i'(z)={m}_i$, where $(1:m_2:\ldots:m_n:0)$ is the infinity
point associated to $B$. In this case $m_i=0$, so there is
$\delta_0>0$ such that for $|z|>\delta_0$, it holds that
$|r_i'(z)|<1/\sqrt{2}$.
 Hence, applying the Mean Value Theorem  (see \cite{Ahlfors}), 
 we have that  if $|x|,|y|>\delta_0$, then
$$|{r}_i(x)-{r}_i(y)|^2=(\Re ({r}_i\,'(c_1))^2+\Im ({r}_i\,'(c_2))^2)|x-y|^2,\quad i=2,\ldots,n$$
where $\Re(q)$ and $\Im(q)$ denote the real part and the imaginary part of $q(z)\in {{\Bbb C}\ll z\gg}$,   respectively,  and $c_1, c_2\in ]x, y[$, where $]x, y[:= \{z \in {\Bbb C}:\, z =x+ (x- y)t,\, t\in (0,1)\}.$
Since $|r_i'(z)|<1/\sqrt{2}$ for $|z|>\delta_0$, we get that $|r_i(y)-r_i(x)|<|x-y|$, for $i=2,\ldots,n$.  In
addition, since $|y|>|x|-K$, we deduce that $|r_i(y)-r_i(x)|<|x-y|$ for
$|x|>\delta_0+K$, and $ i=2,\ldots,n$.\\
\noindent Now, substituting in (\ref{Eq-ultimo-lema}), we get that
$$|r_i(x)-\overline{r}_i(y)|>|\overline{r}_i(y)-r_i(y)|-|x-y|$$
which implies that $\phi(x,y)>|\overline{r}_i(y)-r_i(y)|$ for
$i=2,\ldots,n$. Note that, since $B$ and $\overline{B}$ are divergent
branches, there exists $i_0\in \{1,\ldots,n\}$ such that
$|\overline{r}_{i_0}(y)-r_{i_0}(y)|$ may be as large as we want by
choosing $|x|$ (and thus $|y|$) large enough (see Remark
\ref{R-distance0}, statement 2). Then, for each $K>0$, there exists
$\delta>0$ such that if $|x|>\delta$, it holds that
$\phi(x,y)>|\overline{r}_{i_0}(y)-r_{i_0}(y)|>K.$
\hfill $\Box$

\para
\para

Under these conditions, we obtain Theorem
\ref{P-CharacterizationHausdorff}  that characterizes whether the
Hausdorff distance between two curves is finite.

\para
\para

\begin{theorem}\label{P-CharacterizationHausdorff}
Let ${\cal C}$ and ${\overline{{\cal C}}}$ be two algebraic space curves. It holds that  ${\cal C}$ and ${\overline{{\cal C}}}$   have
a similar asymptotic behavior if and only if the Hausdorff distance between
them is finite.
\end{theorem}

\noindent\textbf{Proof:} First, let us prove that if  ${\cal C}$ and
${\overline{{\cal C}}}$   have a similar asymptotic behavior then,
the Hausdorff distance between them is finite.

\para

\noindent
Let $\kappa$ be the number of infinity branches
of ${\cal C}$. Then,
${\cal C}=B_1\cup \cdots \cup B_{\kappa} \cup \widehat{B},$
where $\widehat{B}$ is the set of points of ${\cal C}$ that do not
belong to any infinity branch. Thus,
$$\sup_{p\in {\cal C}}d(p,{\overline{{\cal C}}})=\max\{\sup_{p\in {B_1}}d(p,{\overline{{\cal C}}}),...,\sup_{p\in {B_{\kappa}}}
d(p,{\overline{{\cal C}}}),\sup_{p\in
\widehat{B}}d(p,{\overline{{\cal C}}})\}.$$ For each
$i=1,...,{\kappa}$, let $B_i=\bigcup_{j=1}^{N_i}L_{i,j}$, where
$L_{i,j}=\{(z,r_{i,j,2}(z),\ldots,r_{i,j,n}(z))\in {\Bbb C}^n:\,z\in {\Bbb C},\,|z|>M_i\}$, and
$N_i=\nu(B_i)$. Then,
$$\sup_{p\in {B_i}}d(p,{\overline{{\cal C}}})=\max_{j=1,\ldots,N_i}\left\{\sup_{|z|>M_i}d((z,r_{i,j,2}(z),\ldots,r_{i,j,n}(z)),{\overline{{\cal C}}})\right\}.$$
Moreover, since ${\cal C}$ and ${\overline{{\cal C}}}$   have a
similar asymptotic behavior then there exists an infinity branch
$\overline{B}_i\subseteq \overline{{\cal C}}$ non-divergent with
$B_i$ (see Definition \ref{D-sim-behavior}). This implies that there
is a leaf
$${\overline L}_{i,j}=\{(z,\overline{r}_{i,j,2}(z),\ldots,\overline{r}_{i,j,n}(z))\in {\Bbb C}^n:\,z\in {\Bbb C},\,|z|>\overline{M}_i\}\subset {\overline{B}_i}$$
such that
$$\lim_{z\rightarrow\infty} d((r_{i,j,2}(z),\ldots,r_{i,j,n}(z)), (\overline{r}_{i,j,2}(z),\ldots,\overline{r}_{i,j,n}(z))=c_{i,j}<\infty$$
(see Lemma \ref{L-non-div-branches} and Remark \ref{R-limit-c}).
Then
$$\lim_{z\rightarrow\infty}d((z,r_{i,j,2}(z),\ldots,r_{i,j,n}(z)),\overline{\cal C})\leq $$$$ \lim_{z\rightarrow\infty}d((z,r_{i,j,2}(z),\ldots,r_{i,j,n}(z)),(z, \overline{r}_{i,j,2}(z),\ldots,\overline{r}_{i,j,n}(z)))=c_{i,j}<\infty$$
Hence, given   $ \eta>0$ there exists $\delta>0$ such that for
$|z|>\delta$ it holds that
$$d((z,r_{i,j,2}(z),\ldots,r_{i,j,n}(z)),\overline{\cal C})<\eta $$
for every $i=1,\ldots,{\kappa}$ and $j=1,\ldots,N_i$.

\para

On the other hand, since $r_{i,j,2},\ldots, r_{i,j,n}$ are
continuous functions, and $\{z\in {\Bbb C}:\,M_i\leq|z|\leq\delta\}$
is a compact set, there exists $\xi>0$ such that
$$\sup_{M_i\leq|z|\leq\delta}d((z,r_{i,j,2}(z),\ldots,r_{i,j,n}(z)),
{\overline{{\cal C}}})<\xi$$ for every $i=1,\ldots,{\kappa}$ and
$j=1,\ldots,N_i$.

\para

\noindent As a consequence, we have that $$\sup_{p\in
{B_i}}d(p,{\overline{{\cal C}}})\leq \max\{\xi,\eta\}<\infty.$$

\noindent
Now, let  $p=(a_1,\ldots, a_n)\in \widehat{B}$. From Lemma
\ref{L-BolaMaxima}, we have that there exists  $K\in \mathbb{R}^+$
such that $|a_i|\leq K$, for $i=1,\ldots,n$. Thus, $d(p,{\cal O})\leq K$, where ${\cal O}$ is the
origin and, $$d(p,{\overline{{\cal C}}})\leq d(p,{\cal O})+d({\cal O},{\overline{{\cal C}}})\leq K+d({\cal O},{\overline{{\cal C}}}).$$
Note that $K<\infty$, and $d({\cal O},{\overline{{\cal C}}})<\infty,$  which implies that $\sup_{p\in \widehat{B}}d(p,{\overline{{\cal C}}})<\infty$.

\para

\noindent Therefore, we conclude that $\sup_{p\in {\cal C}}d(p,{\overline{{\cal C}}})<\infty$. Reasoning similarly, we deduce that $\sup_{\overline{p}\in
{\overline{{\cal C}}}}d(\overline{p},{\cal C})<\infty$, which implies that $d_H({\cal C},\overline{\cal C})<\infty$.

\para
\para
\para

Reciprocally, let us
assume that the Hausdorff distance between ${\cal C}$ and
${\overline{{\cal C}}}$ is finite (that is, $d_H({\cal
C},{\overline{{\cal C}}})=K<\infty$), and let us  prove that  the asymptotic behavior of both curves is similar (i.e.
 for any infinity branch $B\subseteq {\cal C}$  there exists an
infinity branch $\overline{B}\subseteq{\overline{{\cal C}}}$ that
does not diverge with $B$).

\para

For this purpose, we assume that this statement does not hold and let
$B=\{(z,r_{i,j,2}(z),\ldots,r_{i,j,n}(z))\in {\Bbb C}^n:\,z\in {\Bbb
C},\,|z|>M\}\subseteq {\cal C}$ be such that every infinity branch
of ${\overline{{\cal C}}}$ diverges from $B$. Then, according to
Lemma \ref{L-divergent-branches}, for each infinity branch
$\overline{B}_i=\{(z,\overline{r}_{i,j,2}(z),\ldots,\overline{r}_{i,j,n}(z))\in
{\Bbb C}^n:\,z\in {\Bbb C},\,|z|>\overline{M}_i\}\subseteq
{\overline{{\cal C}}}$  $(i=1,\ldots,\kappa)$, there exists
$\delta_i>0$ such that if $|x|>\delta_i$, then
$$d((x,r_{i,j,2}(x),\ldots,r_{i,j,n}(x)),(\overline{a}_1,\overline{a}_2,\ldots,\overline{a}_n))>K$$
for every $(\overline{a}_1,\overline{a}_2,\ldots,\overline{a}_n)\in
\overline{B}_i$. In addition, from Lemma \ref{L-BolaMaxima}, there
exists  $\delta_0>0$ such that any point
$(\overline{a}_1,\overline{a}_2,\ldots,\overline{a}_n)\in{\overline{{\cal
C}}}$ with $|\overline{a}_j|>\delta_0$ for some $j=1,\ldots,n$,
belongs to some infinity branch $\overline{B}_i\subseteq
{\overline{{\cal C}}}$.

\para

\noindent Under these conditions, let
$\delta:=\max\{\delta_0,\delta_1,\ldots,\delta_{\kappa}\}$, and we
consider a point $(x,r_{i,j,2}(x),\ldots,r_{i,j,n}(x))\in B$ such
that $|x|>\delta+K$. Since $d_H({\cal C},{\overline{{\cal C}}})=K$,
there should exist some point
$(\overline{a}_1,\overline{a}_2,\ldots,\overline{a}_n)\in{\overline{{\cal
C}}}$ such that
$$d((x,r_{i,j,2}(x),\ldots,r_{i,j,n}(x)),(\overline{a}_1,\overline{a}_2,\ldots,\overline{a}_n))\leq
K.$$ However, this implies that $|\overline{a}_1|>|x|-K$  (see the
proof of Lemma \ref{L-divergent-branches}) and, hence,
$|\overline{a}_1|>\delta$. Now, Lemma \ref{L-BolaMaxima} states that
this point must belong to some infinity branch
$\overline{B}_i\subseteq{\overline{{\cal C}}}$ and then, Lemma
\ref{L-divergent-branches} claims that
$$d((x,r_{i,j,2}(x),\ldots,r_{i,j,n}(x)),(\overline{a}_1,\overline{a}_2,\ldots,\overline{a}_n))>
K,$$ which is a contradiction. \hfill $\Box$

\para

\para
\para
The following algorithm allows us to decide whether the Hausdorff
distance between two curves $\cal C$ and $\overline{{\cal C}}$ is
finite.  We assume that we have  prepared   $\cal C$ and
$\overline{{\cal C}}$  by means  of a suitable linear change of
coordinates (the same change applied to both curves), such that
$(0:a_2:\ldots: a_n:0)$ ($a_i\not=0$ for some $i=2,\ldots,n$)  is
not an infinity point of ${\cal C}^*$ and $\overline{{\cal C}}^*$
(see Remark \ref{R-infinitypoint}).


\noindent
\begin{center}
\fbox{\hspace*{2 mm}\parbox{13.2cm}{ \vspace*{2 mm} {\bf Algorithm
{\sf Hausdorff Distance.}}
\vspace*{0.2cm}

\noindent Given two   algebraic space curves $\cal C$ and
$\overline{{\cal C}}$  in the $n$-dimensional space, the algorithm
decides whether the Hausdorff distance between  $\cal C$ and $\overline{{\cal C}}$
is finite.
\begin{enumerate}
\item[1.] Compute the infinity points of $\cal C$ and $\overline{{\cal C}}$. If they are
not the same,  {\sc Return} {\it  the Hausdorff distance between the curves $\cal C$ and $\overline{{\cal C}}$ is not finite}. Otherwise, let
$P_1,\ldots,P_{\kappa}$ be these infinity points.
\item[2.]
For each $P_{\ell}:=(1:m_{\ell,2}:\ldots:
m_{\ell,n}:0)$,\,$\ell=1,\ldots, \kappa$ do:
\begin{enumerate}
\item[2.1.] Compute the infinity branches of $\cal C$ associated
to $P_{\ell}$ (see Subsections \ref{Sub-Impli} and \ref{Sub-Param}).
Let $B_1,...,B_{n_{\ell}}$ be these branches. For each
$i=1,\ldots,n_{\ell}$, let $B_i=\{(z,r_{i,2}(z),\ldots,
r_{i,n}(z))\in {\Bbb C}^n:\,z\in {\Bbb C},\,|z|>M_i\}$.
\item[2.2.] Compute the infinity branches  of $\overline{{\cal C}}$ associated
to $P_{\ell}$  (see Subsections \ref{Sub-Impli} and
\ref{Sub-Param}). Let $\overline{B}_1,...,\overline{B}_{l_{\ell}}$
be these branches. For each $j=1,\ldots,l_{\ell}$, let
$\overline{B}_j=\{(z,\overline{r}_{j,2}(z),\ldots,
\overline{r}_{j,n}(z))\in {\Bbb C}^n:\,z\in {\Bbb C},\,|z|>M_j\}$.
\item[2.3.] For each $i=1,\ldots,n_{\ell}$, find
  $j=1,\ldots,l_{\ell}$ such that the
terms with positive exponent in $r_{i,k}(z)$ and
$\overline{r}_{j,k}(z)$ for $k=2, \ldots, n$, are the same up to
conjugation. If there isn't such $j=1,\ldots,l_{\ell}$,  {\sc
Return} {\it  the Hausdorff distance between the curves $\cal C$ and $\overline{{\cal C}}$ is not finite}  (see Lemmas \ref{L-DistVertical} and \ref{L-non-div-branches}, and Theorem \ref{P-CharacterizationHausdorff}).
\item[2.4.] For each  $j=1,\ldots,l_{\ell}$, find
   $i=1,\ldots,n_{\ell}$ such that the
terms with positive exponent in  $r_{i,k}(z)$ and
$\overline{r}_{j,k}(z)$ for $k=2, \ldots, n$, are the same up to
conjugation. If there isn't such $i=1,\ldots,n_{\ell}$,  {\sc
Return} {\it  the Hausdorff distance between the curves $\cal C$ and $\overline{{\cal C}}$ is not finite} (see Lemmas \ref{L-DistVertical} and \ref{L-non-div-branches}, and Theorem \ref{P-CharacterizationHausdorff}).
\end{enumerate}
\item[3.] {\sc Return}  {\it  the Hausdorff distance between the curves $\cal C$ and $\overline{{\cal C}}$ is  finite}.
\end{enumerate}}\hspace{2 mm}}
\end{center}

\para

\para

\para

In the following, we illustrate  the performance of algorithm {\sf
Hausdorff Distance} with two examples. In the first one, we compare
two rational curves defined parametrically. In the second one, the
curves are defined implicitly.

\para

\begin{example}\label{Ej-Haussdorf}
Let $\cal C$  and $\overline{{\cal C}}$ be two rational space curves
in the $4$-dimensional space defined  by the parametrizations
$${\cal P}(s)=\left(\frac{-1+2s^3-s}{s}, \frac{s+1}{s}, \frac{-1}{s}, \frac{s^2+3s-5}{s}\right)$$
and $$\overline{\cal P}(s)=\left(\frac{-1+2s^3-s^2}{s},
\frac{s+1}{s}, \frac{-1}{s}, \frac{s^2+3s-5}{s}\right),$$
respectively. We apply  the algorithm {\sf Hausdorff Distance} to
decide whether the Hausdorff distance between $\cal C$  and
$\overline{{\cal C}}$  is finite.

\begin{itemize}
\item[] \mbox{\sf {Step 1}:}  Compute
the infinity points of  $\cal C$ and $\overline{{\cal C}}$. We
obtain that $\cal C$  and $\overline{{\cal C}}$  have the same
infinity points: $P_1=(1:-1:1:5:0)$ and $P_2=(1:0:0:0:0)$.

\para

We start by analyzing  the infinity branches associated to $P_1$:
\item[]  \mbox{\sf {Step 2.1}:} Reasoning as in Example \ref{E-infbranches-para}, we get only one infinity branch associated to $P_1$
in $\cal C$. It is given by
$B_1=\{(z,r_{1,2}(z),r_{1,3}(z),r_{1,4}(z))\in {\Bbb C}^4:\,z\in
{\Bbb C},\,|z|>M_1\}$, where
$$r_{1,2}(z)=-z+2z^{-2}-4z^{-3}-13z^{-4}-11z^{-5}+\cdots,$$
$$r_{1,3}(z)=z+1-2z^{-2}+4z^{-3}+13z^{-4}+11z^{-5}+\cdots,$$
$$r_{1,4}(z)=5z+8-z^{-1}-9z^{-2}+19z^{-3}+64z^{-4}+62z^{-5}+\cdots.$$

\item[]  \mbox{\sf {Step 2.2}:} We also have that there exists only one infinity branch associated to $P_1$ in $\overline{{\cal C}}$. It  is given by
$\overline{B}_1=\{(z,\overline{r}_{1,2}(z),\overline{r}_{1,3}(z),\overline{r}_{1,4}(z))\in
{\Bbb C}^4:\,z\in {\Bbb C},\,|z|>\overline{M}_1\}$, where
$$\overline{r}_{1,2}(z)=-z+1+z^{-1}+2z^{-2}+z^{-3}-4z^{-4}-7z^{-5}+\cdots,$$
$$\overline{r}_{1,3}(z)=z-z^{-1}-2z^{-2}-z^{-3}+4z^{-4}+7z^{-5}+\cdots,$$
$$\overline{r}_{1,4}(z)=5z+3-6z^{-1}-10z^{-2}-6z^{-3}+18z^{-4}+33z^{-5}+\cdots.$$

\item[]  \mbox{\sf {Step 2.3}} and   \mbox{\sf {Step 2.4}:}
$r_{1,j}(z)$ and $\overline{r}_{1,j}(z),\, j=2,3,4$ have the same terms with positive exponent. Thus, the branches $B_1$ and $\overline{B}_1$ do not diverge.  \\

Now we analyze the infinity branches associated to $P_2$:
\item[]  \mbox{\sf {Step 2.1}:} Reasoning as in Example \ref{E-infbranches-para}, we get that the only infinity branch associated to $P_2$
in $\cal C$ is given by $B_2=\{(z,r_{2,2}(z),r_{2,3}(z),r_{2,4}(z))\in {\Bbb
C}^4:\,z\in {\Bbb C},\,|z|>M_2\}$, where
$$r_{2,2}(z)=1+\sqrt{2}z^{-1/2}-\frac{\sqrt{2}z^{-3/2}}{2}-z^{-2}+\frac{3\sqrt{2}z^{-5/2}}{8}+2z^{-3}+\cdots,$$
$$r_{2,3}(z)=-\sqrt{2}z^{-1/2}+\frac{\sqrt{2}z^{-3/2}}{2}+z^{-2}-\frac{3\sqrt{2}z^{-5/2}}{8}-2z^{-3}+\cdots,$$
$$r_{2,4}(z)=\frac{\sqrt{2}z^{1/2}}{2}+3-\frac{19\sqrt{2}z^{-1/2}}{4}+\frac{z^{-1}}{2}-\frac{39\sqrt{2}z^{-3/2}}{16}+\frac{9z^{-2}}{2}
+\cdots.$$

\noindent We note that $\nu(B_2)=2$, and thus $B_2$  has 2
(conjugated) leaves. That is, $B_2=L_{2,1}\cup L_{2,2}$, where
$L_{2,i}$ are obtained by conjugation in the series $r_{2,2},
r_{2,3}$ and $r_{2,4}$.
\item[]  \mbox{\sf {Step 2.2}:} We also have that there exists only one infinity branch associated to $P_2$ in $\overline{{\cal C}}$. It  is given by
$\overline{B}_2= \{(z,\overline{r}_{2,2}(z),\overline{r}_{2,3}(z),\overline{r}_{2,i4}(z))\in
{\Bbb C}^4:\,z\in {\Bbb C},\,|z|>\overline{M}_2\}$, $i=1,2$, and
$$\overline{r}_{2,2}(z)=1+\sqrt{2}z^{-1/2}-\frac{z^{-1}}{2}-\frac{\sqrt{2}z^{-3/2}}{16}-z^{-2}+\frac{383\sqrt{2}z^{-5/2}}{512}-\frac{z^{-3}}{2}+\cdots,$$
$$\overline{r}_{2,3}(z)=-\sqrt{2}z^{-1/2}+\frac{z^{-1}}{2}+\frac{\sqrt{2}z^{-3/2}}{16}+z^{-2}-\frac{383\sqrt{2}z^{-5/2}}{512}+\frac{z^{-3}}{2}+\cdots,$$
$$\overline{r}_{2,4}(z)=\frac{\sqrt{2}z^{1/2}}{2}+\frac{13}{4}-\frac{159\sqrt{2}z^{-1/2}}{32}+3z^{-1}-\frac{449\sqrt{2}z^{-3/2}}{1024}+5z^{-2}
+\cdots.$$

\noindent
 We note that $\nu(\overline{B}_2)=2$, and thus
$\overline{B}_2$ has 2 (conjugated) leaves. That is,
$\overline{B}_2=\overline{L}_{2,1}\cup \overline{L}_{2,2}$, where
$\overline{L}_{2,i}$ are obtained by conjugation in the series
$\overline{r}_{2,2}, \overline{r}_{2,3}$ and $\overline{r}_{2,4}$.
\item[]  \mbox{\sf  {Step 2.3}} and   \mbox{\sf {Step 2.4}:}
$r_{2,j}(z)$ and $\overline{r}_{2,j}(z),\,j=2,3,4$ have the same
terms with positive exponent. Thus, the branches $B_2$ and
$\overline{B}_2$ do not diverge.

\item[] \mbox{\sf {Step 3}:} The algorithm returns that the Hausdorff distance
between $\cal C$ and $\overline{{\cal C}}$ is finite.\\

\end{itemize}

\para

\para

\para

\noindent We observe that, in this case, the infinity branches of
$\cal C$ and $\overline{{\cal C}}$ do not converge neither diverge
(see Figure \ref{F-Ej-Haussdorf}).

\para

\begin{figure}[h]
$$
\begin{array}{c}
\psfig{figure=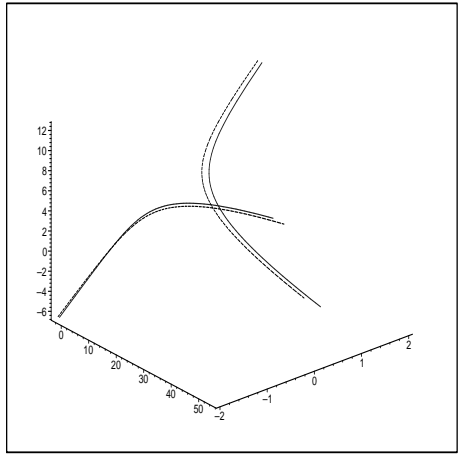,width=6cm,height=6cm,angle=270}
\end{array}
\vspace*{-0.5cm}
$$ \caption{Projections of $\cal C$ and  $\overline{{\cal C}}$ along the axis $x_2$.}\label{F-Ej-Haussdorf}
\end{figure}

\end{example}

\para

\para

\begin{example}\label{Ej-Remark}
Let $\cal C$  and $\overline{{\cal C}}$ be two space curves in the
$3$-dimensional space  implicitly defined  by the polynomials
$$f_1(x_1,x_2,x_3)=-x_2+x_1^2-2x_1x_2^2+x_2^4,\quad f_2(x_1,x_2,x_3)=x_1+x_2^2-x_3x_2^2-x_3$$
and $$\overline{f}_1(x_1,x_2,x_3)=x_2^2-x_1 ,\quad
\overline{f}_2(x_1,x_2,x_3)=2x_1-x_3x_2^2-x_3,$$ respectively. We
apply  the algorithm {\sf Hausdorff Distance} to decide whether the
Hausdorff distance between $\cal C$ and $\overline{{\cal C}}$ is
finite:
\begin{itemize}
\item[] \mbox{\sf {Step 1}:}  Compute
the infinity points of  $\cal C$ and $\overline{{\cal C}}$. We
obtain that $\cal C$  and $\overline{{\cal C}}$  have $P=(1:0:0:0)$
as their unique infinity point.

\para

We analyze  the infinity branches associated to $P$:
\item[]  \mbox{\sf {Step 2.1}:} Reasoning as in Example \ref{E-infbranches-impli}, we get that the only infinity branch associated to $P$
in $\cal C$ is given by $B=\{(z,r_{2}(z),r_{3}(z))\in {\Bbb C}^3:\,z\in {\Bbb
C},\,|z|>M\}$, where
$$r_{2}(z)=z^{1/2}+\frac{z^{-1/4}}{2}-\frac{z^{-7/4}}{64}+\frac{z^{-10/4}}{128}+\cdots,$$
$$r_{3}(z)=2-z^{-3/4}-2z^{-1}+\frac{3z^{-3/2}}{4}+3z^{-7/4}+\cdots.$$
We note that $\nu(B)=4$, and thus $B$  has 4 (conjugated)
leaves. That is, $B=\bigcup_{\alpha=1}^4 L_{\alpha}$, where $L_\alpha$ are obtained by conjugation in the series $r_2$ and $r_3$.
\item[]  \mbox{\sf {Step 2.2}:} We also have that there exists only one infinity branch associated to $P$ in $\overline{{\cal C}}$. It  is given by
$\overline{B}=\{(z,\overline{r}_{2}(z),\overline{r}_{3}(z))\in
{\Bbb C}^3:\,z\in {\Bbb C},\,|z|>\overline{M}\}$, where
$$\overline{r}_{2}(z)=z^{1/2},$$
$$\overline{r}_{3}(z)=2-2z^{-1}+2z^{-2}-2z^{-3}+2z^{-4}-2z^{-5}+\cdots.$$
We note that $\nu(\overline{B})=2$, and thus $\overline{B}$  has 2 (conjugated)
leaves. That is, $\overline{B}=\bigcup_{\beta=1}^2 \overline{L}_{\beta}$, where $\overline{L}_\beta$ are obtained by conjugation in the series $\overline{r}_2$ and $\overline{r}_3$.
\item[]  \mbox{\sf {Step 2.3}} and   \mbox{\sf {Step 2.4}:}
$r_{j}(z)$ and $\overline{r}_{j}(z),\, j=2,3,$ have the same terms
with positive exponent. Thus, the infinity branches $B$ and
$\overline{B}$ do not diverge.
\item[] \mbox{\sf {Step 3}:} The algorithm returns that the Hausdorff distance between the curves $\cal C$
and $\overline{{\cal C}}$ is finite (see Figure \ref{F-Ej-Remark}).
\end{itemize}
 \begin{figure}[h]
$$
\begin{array}{lcr}
\psfig{figure=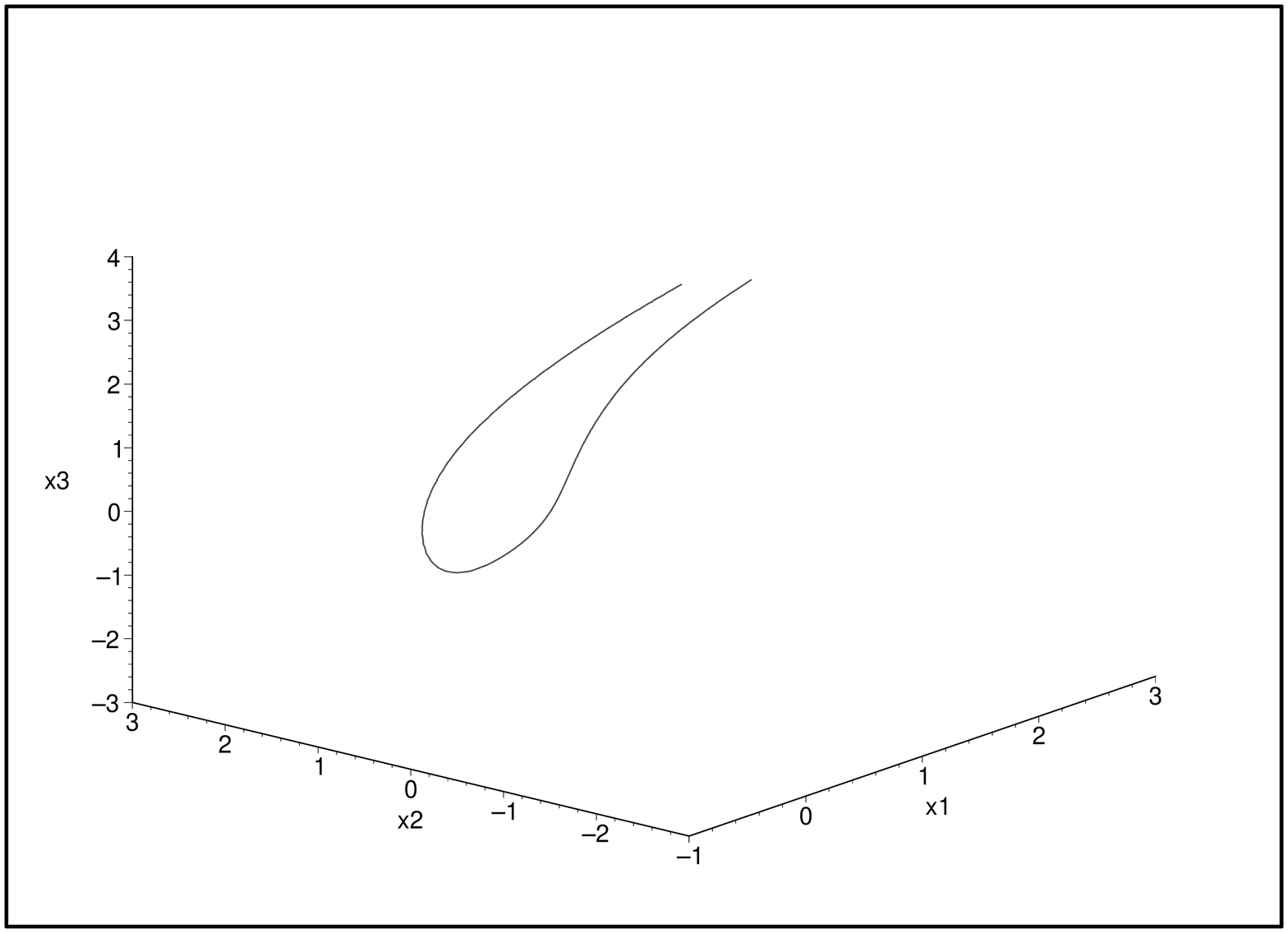,width=4.3cm,height=4.3cm,angle=270}
&
\psfig{figure=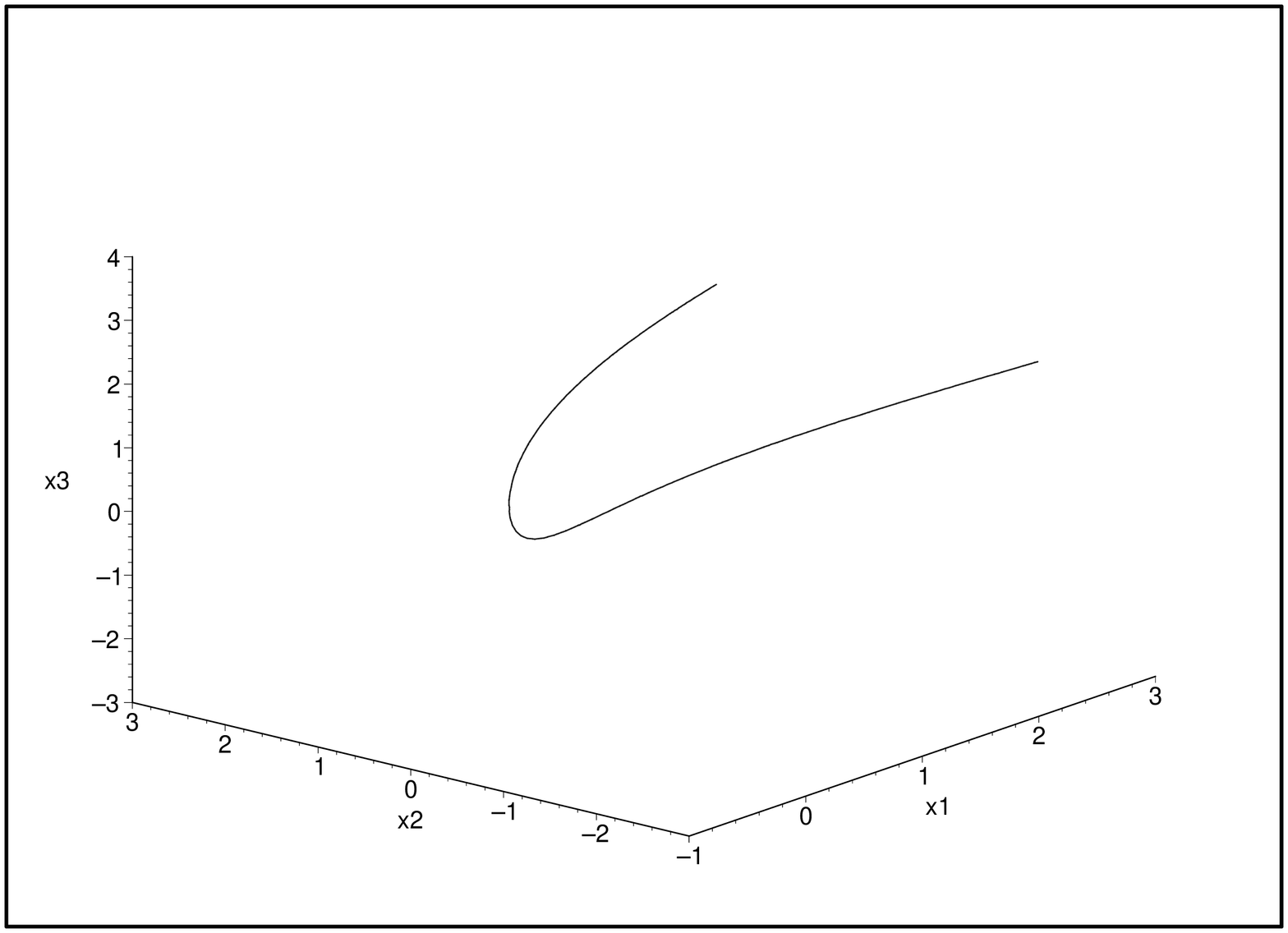,width=4.3cm,height=4.3cm,angle=270}
&
\psfig{figure=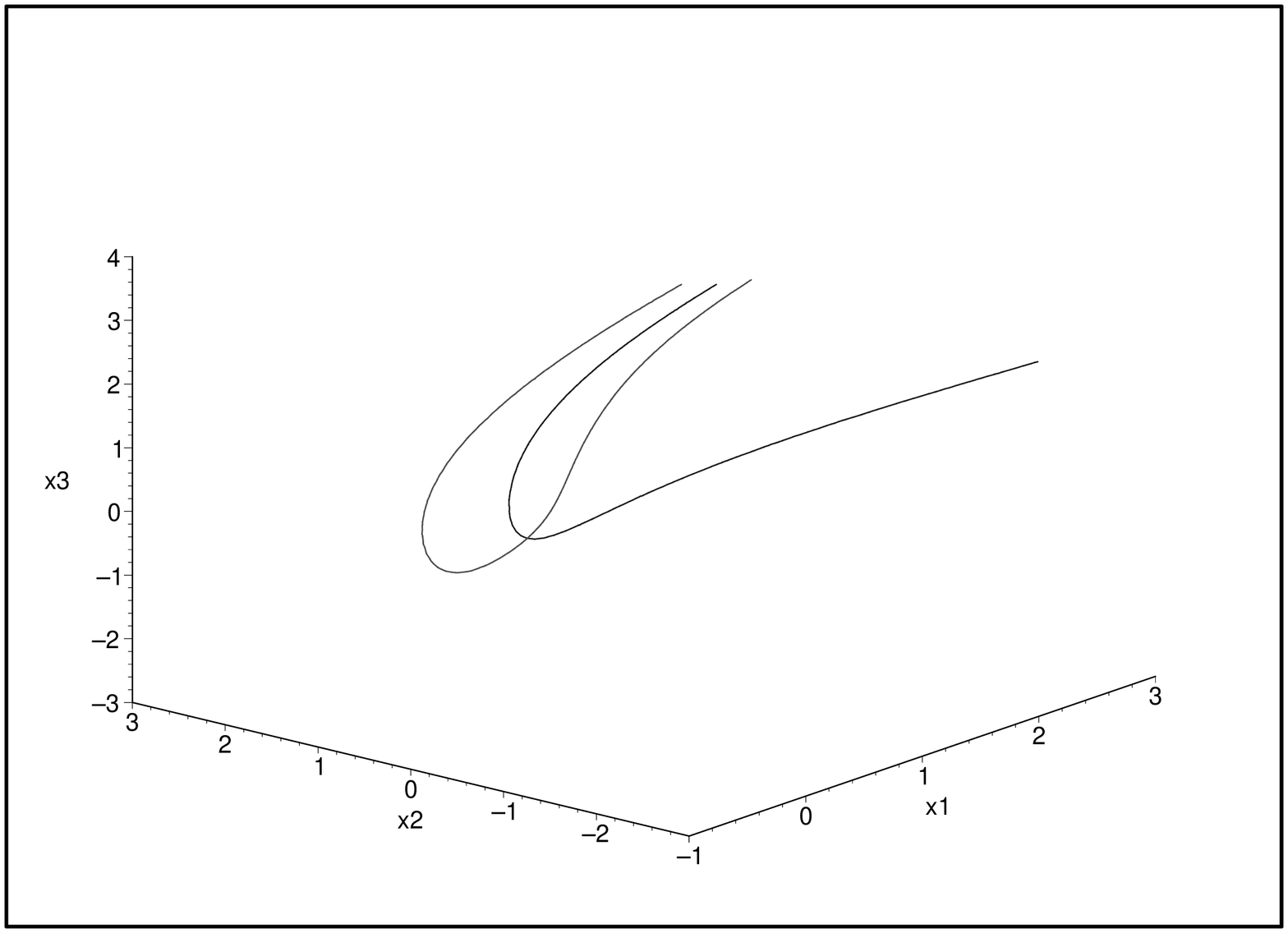,width=4.3cm,height=4.3cm,angle=270}\end{array}
\vspace*{-0.5cm}
$$ \caption{$\cal C$ (left),  $\overline{{\cal C}}$ (center), and the asymptotic behavior of $\cal C$ and $\overline{{\cal C}}$ (right)}
\label{F-Ej-Remark}
\end{figure}
\end{example}

\para

\para
\begin{remark}
Results obtained in Example \ref{Ej-Remark} could be surprising  for
the reader since in Figure \ref{F-Ej-Remark}, the Hausdorff distance
between $\cal C$ and $\overline{{\cal C}}$ does not seem to be
finite. The explanation of this phenomenon is  that throughout this
paper, we are dealing with the whole curve (including its complex
part) but clearly, if we restrict to the real part, the Hausdorff
distance could go from being finite (if we consider the curves over
$\Bbb C$) to be infinite (if we consider the curves over $\Bbb R$).
In this example, the Hausdorff distance is infinity if we restrict
to the real part of the curves. More precisely, in Example
\ref{Ej-Remark},   the infinity branch $B\subseteq{\cal C}$ has two
complex leaves that cannot be plotted. They are
${L}_{3}=\{(z,{r}_{3,2}(z),{r}_{3,3}(z))\in {\Bbb C}^3:\,z\in {\Bbb
C},\,|z|>{M}\}$, where
$$r_{3,2}(z)=-z^{1/2}+\frac{Iz^{-1/4}}{2}+\frac{Iz^{-7/4}}{64}-\frac{z^{-5/2}}{128}+\cdots,$$
$$r_{3,3}(z)=2+Iz^{-3/4}-2z^{-1}-\frac{3z^{-3/2}}{4}-3Iz^{-7/4}+\cdots,$$
and ${L}_{4}=\{(z,{r}_{4,2}(z),{r}_{4,3}(z))\in {\Bbb C}^3:\,z\in
{\Bbb C},\,|z|>{M}\}$, where
$$r_{4,2}(z)=-z^{1/2}-\frac{Iz^{-1/4}}{2}-\frac{Iz^{-7/4}}{64}-\frac{z^{-5/2}}{128}+\cdots,$$
$$r_{4,3}(z)=2-Iz^{-3/4}-2z^{-1}-\frac{3z^{-3/2}}{4}+3Iz^{-7/4}+\cdots.$$
Note that the imaginary parts of these series are given by
terms with negative exponent, which means that they vanish as $z$
grows to infinity. Hence, both leaves converge to the
real leaf
$\overline{L}_{2}=\{(z,\overline{r}_{2,2}(z),\overline{r}_{2,3}(z))\in
{\Bbb C}^3:\,z\in {\Bbb C},\,|z|>\overline{M}\}$, where
$$\overline{r}_{2,2}(z)=-z^{1/2},$$
$$\overline{r}_{2,3}(z)=2-2z^{-1}+2z^{-2}-2z^{-3}+2z^{-4}-2z^{-5}+\cdots,$$
that belongs to the branch $\overline{B}\subseteq\overline{{\cal
C}}$. \para

\noindent Summarizing, Example \ref{Ej-Remark} shows that a complex
leaf may converge to a real one. Furthermore,  the Hausdorff
distance between two curves may be finite while the distance between
their real parts is infinite.
\end{remark}

\end{document}